\documentclass[preprint]{elsarticle}

\usepackage[english]{babel}
\usepackage{hyphenat}
\usepackage{acro}

\usepackage[a4paper, total={6.5in, 8.66in}]{geometry}
\usepackage{tocbibind}
\usepackage[export]{adjustbox}
\usepackage[framemethod=tikz]{mdframed}
\usepackage{float}

\usepackage{longtable}
\usepackage{booktabs}
\usepackage{multirow}
\usepackage{multicol}

\usepackage{amsmath}
\usepackage{amssymb}
\usepackage{mathtools}
\usepackage{amsthm}
\usepackage[per-mode=symbol,
            inter-unit-product=\cdot,
            exponent-product=\cdot,
            separate-uncertainty=true,
            print-unity-mantissa=false,
            table-alignment-mode=none,
            table-alignment=right]{siunitx}
\usepackage{physics}
\usepackage{cases}

\usepackage{algorithm}
\usepackage{algpseudocode}

\usepackage[hypertexnames=false,
            colorlinks=true,
            pdfborderstyle={/S/U/W 0},
            citecolor=blue]{hyperref}
\usepackage[capitalise, noabbrev]{cleveref}

\usepackage{graphicx}
\usepackage{subcaption}

\usepackage[colorinlistoftodos,
            bordercolor=white,
            textsize=footnotesize]{todonotes}

\newfloat{algorithm}{t}{lop}

\everymath{\displaystyle}
\allowdisplaybreaks

\graphicspath{{fig}{standalone}}

\DeclareSIUnit\mmhg{mmHg}

\theoremstyle{definition}

\newtheorem{definition}{Definition}

\DeclareMathOperator{\atantwo}{arctan_2}

\renewcommand{\epsilon}{\varepsilon}

\renewcommand{\tilde}{\widetilde}

\newcommand{\lifex}{\texttt{life\textsuperscript{x}}}
\newcommand{\dealii}{\texttt{deal.II}}

\newcommand{\domain}{\Omega}

\newcommand{\normal}{\mathbf{n}}

\newcommand{\potential}{v}
\newcommand{\ionicvars}{\mathbf{w}}
\newcommand{\actvars}{\mathbf{s}}
\newcommand{\displacement}{\mathbf{d}}
\newcommand{\circvars}{\mathbf{c}}

\newcommand{\fibers}{\mathbf{f}_0}
\newcommand{\sheets}{\mathbf{s}_0}
\newcommand{\normals}{\mathbf{n}_0}
\newcommand{\difftensor}{\mathbf{D}_\mathrm{m}}
\newcommand{\sigmaf}{\sigma_\text{f}}
\newcommand{\sigmas}{\sigma_\text{s}}
\newcommand{\sigman}{\sigma_\text{n}}
\newcommand{\Iion}{I_\text{ion}}
\newcommand{\Iapp}{I_\text{app}}
\newcommand{\calcium}{[\text{Ca}^{2+}]_\text{i}}

\newcommand{\Itensor}{\mathbf{I}}
\newcommand{\Ftensor}{\mathbf{F}}
\newcommand{\jacobian}{J}
\newcommand{\SL}{\text{SL}}
\newcommand{\actstress}{\mathbf{P}_\text{act}}

\newcommand{\stress}{\mathbf{P}}
\newcommand{\Kepi}{\mathbf{K}_\text{epi}}
\newcommand{\Cepi}{\mathbf{C}_\text{epi}}

\newcommand{\src}{^\text{src}}
\newcommand{\dst}{^\text{dst}}
\newcommand{\srcpoint}[1]{\mathbf x\src_{#1}}
\newcommand{\dstpoint}[1]{\mathbf x\dst_{#1}}
\newcommand{\rbf}{\phi}
\newcommand{\rbfradius}[1]{r_{#1}}
\newcommand{\rbfradiusmax}{r_{\max}}
\newcommand{\interpmatrix}{\Phi_\text{int}}
\newcommand{\evalmatrix}{\Phi_\text{eval}}

\newcommand{\mesh}{\mathcal{M}}

\newcommand{\interp}[1]{\Pi_{#1}}
\newcommand{\interpresc}[1]{\Pi^\text{res}_{#1}}

\newcommand{\geodist}{g}
\newcommand{\discgeodist}{g_h}
\newcommand{\physicaldim}{n}
\newcommand{\realspace}{\mathbb{R}^\physicaldim}
\newcommand{\errinfty}{e_{\infty}}

\newcommand{\rlrbfg}{RL-RBF-G}

\newcommand{\ttp}{_\text{TTP}}
\newcommand{\crn}{_\text{CRN}}

\newcommand{\ra}{_\text{RA}}
\newcommand{\rv}{_\text{RV}}
\newcommand{\la}{_\text{LA}}
\newcommand{\lv}{_\text{LV}}
\newcommand{\ao}{_\text{Ao}}
\newcommand{\pt}{_\text{PT}}
\newcommand{\ventricles}{_\text{V}}
\newcommand{\myo}{_\text{myo}}

\definecolor{plot1}{HTML}{1f77b4}
\definecolor{plot2}{HTML}{ff7f0e}
\definecolor{plot3}{HTML}{2ca02c}
\definecolor{plot4}{HTML}{d62728}
\definecolor{plot5}{HTML}{9467bd}

\acsetup{make-links  = false}

\DeclareAcronym{RBF}{long={radial basis function}, short={RBF}}
\DeclareAcronym{RLRBF}{long={rescaled localized radial basis function}, short={RL-RBF}}
\DeclareAcronym{ODE}{long={ordinary differential equation}, short={ODE}}
\DeclareAcronym{SVD}{long={singular value decomposition}, short={SVD}}
\DeclareAcronym{IMEX}{long={implicit-explicit}, short={IMEX}}
\DeclareAcronym{DoF}{long={degree of freedom}, short={DoF}, long-plural-form={degrees of freedom}}
\DeclareAcronym{ILU}{long={incomplete LU}, short={ILU}}
\DeclareAcronym{TTP}{long={ten Tusscher-Panfilov}, short={TTP06}}
\DeclareAcronym{CRN}{long={Courtemanche-Ramirez-Nattel}, short={CRN}}
\DeclareAcronym{MEF}{long={mechano-electric feedback}, short={MEF}}
\DeclareAcronym{LDRBM}{long={Laplace-Dirichlet rule based method}, short={LDRBM}}

\makeatletter
\def\ps@pprintTitle{%
 \let\@oddhead\@empty
 \let\@evenhead\@empty
 \def\@oddfoot{}%
 \let\@evenfoot\@oddfoot}
\makeatother
\biboptions{sort&compress}

\begin{document}

\begin{frontmatter}
  \title{Robust radial basis function interpolation based on geodesic distance for the numerical coupling of multiphysics problems}

  \author[1]{Michele Bucelli\texorpdfstring{\corref{cor1}}{}}
  \author[1]{Francesco Regazzoni}
  \author[1]{Luca Dede'}
  \author[1,2]{Alfio Quarteroni}

  \affiliation[1]{
    organization={MOX Laboratory of Modeling and Scientific Computing, Dipartimento di Matematica, Politecnico di Milano},
    addressline={Piazza Leonardo da Vinci 32},
    postcode={20133},
    city={Milano},
    country={Italy}}
  \affiliation[2]{
    organization={Institute of Mathematics, École Polytechnique Fédérale de Lausanne},
    addressline={Station 8, Av. Piccard},
    postcode={CH-1015},
    city={Lausanne},
    country={Switzerland (Professor Emeritus)}}

  \cortext[cor1]{Corresponding author. E-mail: michele.bucelli@polimi.it}
  \date{Last update: {\today}}

  \journal{}

  \begin{abstract}
    Multiphysics simulations frequently require transferring solution fields between subproblems with non-matching spatial discretizations, typically using interpolation techniques. Standard methods are usually based on measuring the closeness between points by means of the Euclidean distance, which does not account for curvature, cuts, cavities or other non-trivial geometrical or topological features of the domain. This may lead to spurious oscillations in the interpolant in proximity to these features. To overcome this issue, we propose a modification to \ac{RLRBF} interpolation to account for the geometry of the interpolation domain, by yielding conformity and fidelity to geometrical and topological fieatures. The proposed method, referred to as \rlrbfg{}, relies on measuring the geodesic distance between data points. \rlrbfg{} removes all oscillations appearing in the \ac{RLRBF} interpolant, resulting in increased accuracy in domains with complex geometries. We demonstrate the effectiveness of \rlrbfg{} interpolation through a convergence study in an idealized setting. Furthermore, we discuss the algorithmic aspects and the implementation of \rlrbfg{} interpolation in a distributed-memory parallel framework, and present the results of a strong scalability test yielding nearly ideal results. Finally, we show the effectiveness of \rlrbfg{} interpolation in multiphysics simulations by considering an application to a whole-heart cardiac electromecanics model.
  \end{abstract}

\end{frontmatter}

\acresetall


\section{Introduction}

We consider the problem of transferring data between two different meshes of the same spatial domain. This problem often arises in the context of finite elements multiphysics simulations \cite{regazzoni2022cardiac,piersanti20223d,zappon2022reduced,zappon2023efficient}. Indeed, two or more coupled problems defined on the same domain may have different accuracy requirements, which should be addressed with different, possibly non-matching discretizations, to balance computational cost and accuracy \cite{farah2016volumetric}, or to allow the coupling of independent solvers \cite{chourdakis2022precice}. In all these cases, it is therefore crucial to introduce suitable intergrid transfer operators to realize the coupling. Operators of this kind are relevant also for remeshing techniques or geometric multigrid methods \cite{kronbichler2019multigrid, clevenger2020flexible}.

\begin{figure}
  \centering

  \begin{subfigure}[b]{0.49\textwidth}
    \includegraphics[width=\textwidth]{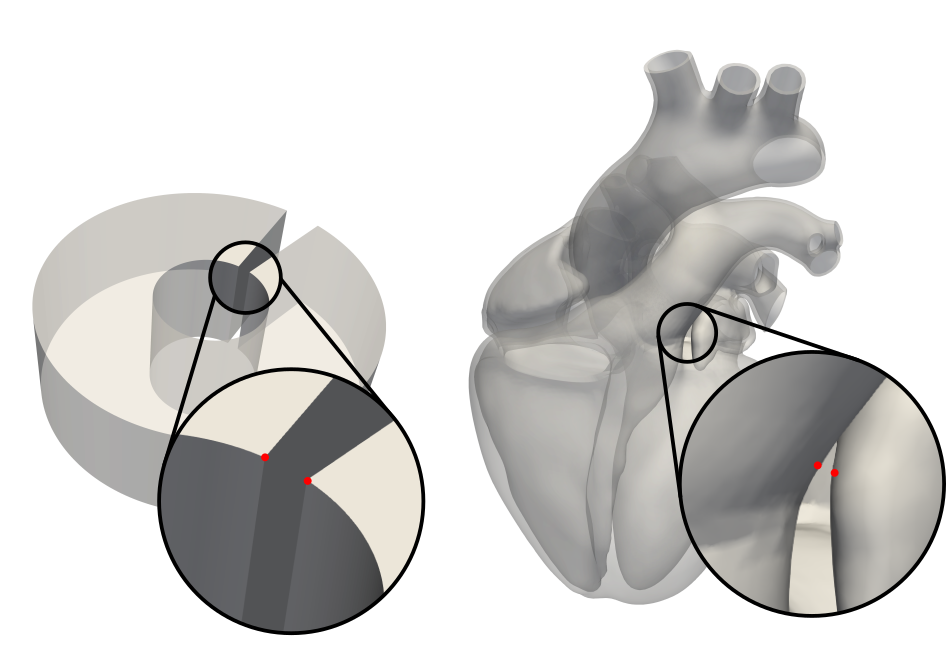}
    \caption{}
    \label{fig:complex-domains}
  \end{subfigure}
  \begin{subfigure}[b]{0.245\textwidth}
    \centering
    \includegraphics{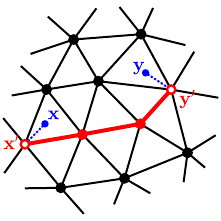}
    \caption{}
    \label{fig:geodesic-distance}
  \end{subfigure}
  \begin{subfigure}[b]{0.245\textwidth}
    \centering
    \includegraphics{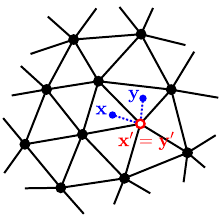}
    \caption{}
    \label{fig:geodesic-distance-positivity}
  \end{subfigure}

  \caption{(a) Examples of domains with complex geometrical features. The red dots mark points on the domains that are close in the Euclidean metric, but that belong to topologically distant regions. The heart domain on the right is the Zygote heart model \cite{zygote}. (b) A schematic representation of the discrete geodesic distance. (c) An example where, even if $\mathbf x \neq \mathbf y$, $\discgeodist(\mathbf x, \mathbf y) = 0$.}
\end{figure}

To this purpose, we consider \ac{RLRBF} interpolation, introduced in \cite{deparis2014rescaled} and previously employed for interface-coupled \cite{voet2022internodes,deparis2016internodes} and volume-coupled \cite{salvador2020intergrid,bucelli2023preserving} problems. When constructing the \ac{RLRBF} interpolant, the closeness between points is measured using the Euclidean norm. This may prevent from capturing the correct geometry of the domain under consideration. Indeed, in realistic applications, disconnected regions of the domain may be close in the Euclidean norm (see e.g. \cref{fig:complex-domains}). Therefore, an interpolant constructed by relying exclusively on the Euclidean distance will not account for the fact that those regions are disconnected. When considering data representing the solution of a physical problem, the points in regions of the domain that are not connected may feature very different solution values, which can lead to large artificial oscillations in the interpolant. These oscillations mostly appear in proximity to non-trivial geometrical features such as holes or cuts in the domain, as showcased by the examples of this paper.

We propose a modification to \ac{RLRBF} interpolation that relies on an approximation of the geodesic distance within the domain of interest to take its geometry into account. The geodesic distance is used to threshold the Euclidean distance, so that points that belong to disconnected regions of the domain are regarded as infinitely distant. We discuss the computational aspects of the proposed method by describing in detail the algorithms used to construct and evaluate the interpolant. We focus in particular on a distributed-memory parallel implementation, in which the interpolation and evaluation points are distributed across several processes.

Finally, we present numerical examples that illustrate the issues of \ac{RLRBF} interpolation in topologically complex domains, and we show that our geodesic distance thresholding method can overcome those issues. We consider both an idealized example and a realistic multiphysics application to three-dimensional whole-heart electromechanics \cite{fedele2023comprehensive, gerach2021electro, augustin2016patient, augustin2021computationally, gurev2011models, washio2013multiscale, levrero2020sensitivity, hirschvogel2017monolithic, karabelas2022accurate, strocchi2020simulating, viola2023gpu, gerach2024differential}. The numerical examples are used to assess the computational efficiency and parallel scalability of the proposed method.

This paper is structured as follows. In \cref{sec:geometry-aware-rbf} we recall the \ac{RLRBF} intergrid transfer and describe the proposed modification leading to the \rlrbfg{} method. In \cref{sec:algorithms} we discuss the computational and algorithmic details of our implementation. \Cref{sec:numerical-experiments} describes numerical experiments that showcase the properties of the proposed method. Finally, we draw some conclusive remarks in \cref{sec:conclusions}.

\section{RBF intergrid transfer with geodesic thresholding}
\label{sec:geometry-aware-rbf}

We briefly recall the \ac{RLRBF} interpolation procedure, and present our novel geodesic-based distance thresholding method to interpolate in complex domains accounting for their geometry.

\subsection{Intergrid transfer using rescaled, localized RBF interpolation}

Let $\Omega \subset \realspace$ be an open, bounded domain, and let $\mesh\src$ and $\mesh\dst$ be two meshes of $\Omega$. The two meshes are not necessarily conforming or nested, and may feature different resolutions and different element types (tetrahedral or hexahedral). In other words, the two meshes represent two different approximations of the domain $\Omega$, which we will denote by $\Omega\src$ and $\Omega\dst$ respectively. We remark that, in general, the boundaries of $\Omega\src$ and $\Omega\dst$ do not necessarily coincide.

Let $f:\Omega\src \to \mathbb{R}$ be a known function. Let $\{\srcpoint{i}\}_{i = 1}^{N\src}$, with $\srcpoint{i} \in \Omega\src$, be a set of interpolation points, and $f\src_i = f(\srcpoint{i})$ for $i = 1, \dots, N\src$. The \ac{RBF} interpolant of the function $f$ at points $\srcpoint{i}$ is a function $\interp{f}:\realspace \to \mathbb{R}$ in the form
\begin{equation}
  \interp{f}(\mathbf x) = \sum_{j = 1}^{N\src} \gamma_j\,\rbf(d(\mathbf x, \srcpoint{j}), \rbfradius{j})\;,
  \label{eq:interpolant}
\end{equation}
wherein $\gamma_j \in \mathbb{R}$ are unknown interpolation coefficients, $d:\realspace \times \realspace \to [0, +\infty)$ measures the distance between two points, and $\rbf:\mathbb{R} \times [0, +\infty) \to [0, 1]$ is the compactly supported $C^2$-continuous Wendland basis function \ac{RBF} \cite{wendland1995piecewise} with support radius $\rbfradius{j} > 0$, defined as
\begin{equation}
  \rbf(t, r) = \max\left\{1 - \frac{t}{r}, 0\right\}^4 \left(1 + 4\frac{t}{r}\right)\quad \text{for } r \geq 0\;.
\end{equation}
The distance is typically measured in the Euclidean norm \cite{deparis2014rescaled, voet2022internodes, wendland1999meshless, morse2005interpolating, franke1998solving, chourdakis2022precice}, that is
\begin{equation}
  d(\mathbf x, \mathbf y) = \|\mathbf x - \mathbf y\|\;.
  \label{eq:distance-euclidean}
\end{equation}
However, as discussed in the following, this may introduce artificial oscillations when the domain features a complex geometry.

The coefficients $\gamma_j$ in \eqref{eq:interpolant} are determined by imposing the interpolation constraints $\interp{f}(\srcpoint{i}) = f\src_i$ for all $i = 1, \dots, N\src$, which is equivalent to solving the following linear system:
\begin{equation}
  \interpmatrix \boldsymbol\gamma = \mathbf f\src\;,
  \label{eq:rbf-linear-system}
\end{equation}
wherein $\boldsymbol\gamma = (\gamma_1, \gamma_2, \dots, \gamma_{N\src})^T$ is the vector of interpolation coefficients, the vector $\mathbf f\src = (f\src_1,\allowbreak f\src_2,\allowbreak \dots,\allowbreak f\src_{N\src})^T$ collects the evaluations of $f$ on the interpolation points, and $\interpmatrix$ is an $N\src \times N\src$ matrix of entries
\begin{equation}
  (\interpmatrix)_{ij} = \rbf(d(\srcpoint{i}, \srcpoint{j}), \rbfradius{j})\;.
\end{equation}

The evaluation of the interpolant $\interp{f}$ onto a set of points $\{\dstpoint{i}\}_{i = 1}^{N\dst}$, with $\dstpoint{i} \in \Omega\dst$, can be obtained by computing the matrix-vector product $\mathbf f\dst = \evalmatrix \boldsymbol\gamma$, where $\mathbf f\dst = (\interp{f}(\dstpoint{1}), \interp{f}(\dstpoint{2}), \dots, \interp{f}(\dstpoint{N\dst}))^T$ and $\evalmatrix$ is an $N\dst \times N\src$ matrix of entries
\begin{equation}
  (\evalmatrix)_{ij} = \phi(d(\dstpoint{i}, \srcpoint{j}), \rbfradius{j})\;.
\end{equation}

The \ac{RBF} interpolant may yield oscillations even if the interpolated data is constant \cite{deparis2014rescaled}. To avoid this, we adopt the \ac{RLRBF} interpolation introduced in \cite{deparis2014rescaled}. The \ac{RBF} support radius $\rbfradius{j}$ associated with interpolation node $\srcpoint{j}$ is chosen as the distance between $\srcpoint{j}$ and the $M$-th nearest neighbor in the set $\{\srcpoint{i}\}_{i = 1}^{N\src}$, magnified by a factor $\alpha$ \cite{bucelli2023preserving,voet2022internodes}, where $M \in \mathbb{N}$ and $\alpha \in [1, \infty)$ are user-defined parameters. Furthermore, introducing the constant function $g(\mathbf x) = 1$, we rescale $\interp{f}$ by the interpolant $\interp{g}$ of $g$, leading to the following definition of the rescaled interpolant $\interpresc{f}$:
\begin{equation}
  \interpresc{f}(\mathbf x) = \frac{\interp{f}(\mathbf x)}{\interp{g}(\mathbf x)}\;,
\end{equation}
To keep the notation light, we shall drop the superscript ``res'' from here on, and denote the rescaled interpolant simply by $\interp{f}$.

The definition of the interpolant and its evaluation generalize to vectorial functions by interpolating every component independently. We refer to \cite{deparis2014rescaled,bucelli2023preserving} for further details on \ac{RLRBF} interpolation, and to \cite{demarchi2020convergence} for a discussion of its convergence properties.

\subsection{Geodesic distance thresholding}
\label{sec:geodesic-distance-thresholding}

The \ac{RBF} interpolant is a linear combination of the interpolation coefficients $\gamma_j$, where the weights of the combination are larger the closer the evaluation point $\mathbf x$ is to the interpolation point $\srcpoint{j}$. In qualitative terms, this means that the value of $\interp{f}(\mathbf x)$ is influenced by the value of $f$ at the interpolation points that are close to $\mathbf x$.

This may not be appropriate for complex domains that are characterized by disconnected regions, high curvature or cavities, such as those represented in \cref{fig:complex-domains}. In those cases, the Euclidean distance may not be an adequate measure of the closeness between two points of $\Omega$. This may lead to large artificial oscillations in the interpolant.

To avoid this issue, we need to replace the Euclidean distance \eqref{eq:distance-euclidean} with a notion of distance that better accounts for the geometry of $\Omega$. As a starting point, we consider the geodesic distance within $\Omega$.

\begin{definition}[Geodesic distance]
  Given two points $\mathbf x, \mathbf y \in \Omega$, their geodesic distance $\geodist(\mathbf x, \mathbf y)$ is the length of the shortest curve contained in $\Omega$ and featuring $\mathbf x$ and $\mathbf y$ as its endpoints.
\end{definition}

Since we consider discrete approximations of the domain $\Omega$, we replace the geodesic distance defined above with a discrete counterpart. We begin by defining the distance between mesh vertices.

\begin{definition}[Discrete geodesic distance]
  \label{def:geodist-discrete}
  Let $\mesh$ be a mesh representing an approximation $\tilde{\Omega}$ of $\Omega$. Given two points $\mathbf x, \mathbf y \in \tilde{\Omega}$, vertices of $\mesh$, we define their discrete geodesic distance $\discgeodist(\mathbf x, \mathbf y)$ as the length of the shortest path from $\mathbf x$ to $\mathbf y$ traveling along edges and element diagonals of $\mesh$.
\end{definition}

Applications of the interpolation procedure often require the interpolation or evaluation points to be internal to mesh elements (such as Gauss-Legendre quadrature points \cite{bucelli2023preserving}). Furthermore, the evaluation of the interpolant requires us to compute the distance between interpolation points in $\Omega\src$ and evaluation points in $\Omega\dst$. Therefore, we extend the above definition to arbitrary points in $\realspace$ as follows.

\begin{definition}[Discrete geodesic distance]
  \label{def:geodist-discrete-extended}
  Let $\mesh$ be a reference mesh approximating $\Omega$. Given two points $\mathbf x, \mathbf y \in \realspace$, we set $\discgeodist(\mathbf x, \mathbf y) = \discgeodist(\mathbf x', \mathbf y')$, where $\mathbf x'$ and $\mathbf y'$ are the vertices of $\mesh$ closest to $\mathbf x$ and $\mathbf y$, respectively.
\end{definition}

\Cref{fig:geodesic-distance} displays a schematic representation of the above notions. The computation of the discrete geodesic distance $\discgeodist$ is computationally efficient and easy to implement \cite{dijkstra1959note,lanthier2001approximating} (see also \cref{sec:algorithm-geodist}). However, $\discgeodist$ may be a poor approximation of the exact geodesic distance $\geodist$ \cite{lanthier1997approximating,lanthier2001approximating}. Furthermore, as shown in \cref{fig:geodesic-distance-positivity}, it may happen that $\discgeodist(\mathbf x, \mathbf y) = 0$ even if $\mathbf x \neq \mathbf y$, which may hinder the well-posedness of the interpolation problem \cite{voet2022internodes}. Instead of resorting to more accurate yet more complex (and computationally demanding) discretization methods, we use $\discgeodist$ to apply a threshold to the Euclidean distance, replacing \eqref{eq:distance-euclidean} with the following.

\begin{definition}[Thresholded Euclidean distance]
  Let $\mesh$ be a reference mesh approximating $\Omega$. Given two points $\mathbf x, \mathbf y \in \realspace$, we define
  \begin{equation}
    d(\mathbf x, \mathbf y; R) = \begin{cases}
      +\infty & \text{if } \discgeodist(\mathbf x, \mathbf y) > R\;, \\
      \discgeodist(\mathbf x, \mathbf y) & \text{if } \beta h_\text{max} + \|\mathbf x - \mathbf y\| < \discgeodist(\mathbf x, \mathbf y) \leq R\;, \\
      \|\mathbf x - \mathbf y\| & \text{if } \discgeodist(\mathbf x, \mathbf y) \leq \beta h_\text{max} + \|\mathbf x - \mathbf y\|\;, \\
    \end{cases}
    \label{eq:distance-cutoff}
  \end{equation}
  where the discrete geodesic distance $\discgeodist(\mathbf x, \mathbf y)$ is evaluated with $\mesh$ as a reference mesh, $h_\text{max}$ is the maximum diameter of the elements of $\mesh$, and $\beta > 0$ is a constant coefficient.
\end{definition}

The coefficient $\beta$ has the role of detecting regions of high curvature, where the discrete geodesic distance between points is significantly larger than the Euclidean distance. In those regions it is preferrable to compute $d(\mathbf x, \mathbf y)$ through $\discgeodist(\mathbf x, \mathbf y)$, instead of the Euclidean distance $\|\mathbf x - \mathbf y\|$. Indeed, although $\discgeodist$ is a potentially inaccurate approximation of $g$, this approximation will generally be better than the Euclidean distance. Conversely, in regions of low curvature, the Euclidean distance enjoys better smoothness properties than $\discgeodist$, and is therefore a better choice for evaluating $d(\mathbf x, \mathbf y)$.

Finally, we replace the definition of the interpolant in \eqref{eq:interpolant} with
\begin{equation}
  \interp{f}(\mathbf x) = \sum_{j = 1}^{N\src} \gamma_j\,\rbf(d(\mathbf x, \srcpoint{j}; \rbfradius{j}), \rbfradius{j})\;,
\end{equation}
where the threshold $R$ of \eqref{eq:distance-cutoff} is chosen equal to the \ac{RBF} support radius. The entries of $\interpmatrix$ and $\evalmatrix$ are updated accordingly. To evaluate the discrete geodesic distance, we use the finest of the two meshes $\mesh\src$ and $\mesh\dst$, to obtain the best approximation. We refer to the resulting method as \rlrbfg{} interpolation.

\section{Algorithms and parallel implementation}
\label{sec:algorithms}

We describe in this section the algorithms and data structures implemented to construct and evaluate the \rlrbfg{} interpolant, discussing computational cost and parallel scalability. Our implementation is designed for a distributed-memory parallel architecture, based on the MPI standard.

All numerical methods presented here are available as open source within the core module of \lifex{} \cite{lifex,africa2022lifexcore}, a library based on \dealii{} \cite{dealii,arndt2020dealii,arndt2022dealii9.4} for finite element simulations of the cardiac function \cite{africa2024lifexcfd,africa2023lifexfiber,africa2023lifexep}.

The algorithm presented in this section significantly improves upon the one we previously discussed in \cite{bucelli2023preserving}. The major advancements are the use of R-trees to perform spatial queries and a significant reduction in inter-process communication. We report in \cref{alg:rbf-matrix-assembly} an overview of the proposed algorithm. The rest of this section describes in more detail its individual steps.

\subsection{Coarse representations of local points}

Since we work in a distributed-memory parallel framework, each process $r$ only stores a subset $X\src_r$ of the source points $\{\srcpoint{i}\}_{i = 1}^{N\src}$. However, when building the interpolation matrix $\interpmatrix$, each process needs to evaluate the distance of each of the points it owns from all other points, including some that are owned by other processes. Those points must therefore be retrieved through inter-process communication.

To do so, we build for every process $r$ a coarse representation $B_r \subset \realspace$ such that $X\src_r \subset B_r$. We define $B_r$ as the union of multiple axis-aligned bounding boxes, obtained by successively splitting in half the axis-aligned box containing all the points. Then, the coarse representation is communicated to all other processes. Finally, each process $r$ uses the coarse representation of every other process $p$ to determine which of the points $X\src_r$ must be communicated to $p$. This communication is obtained by means of a point-to-point communication algorithm as implemented by \dealii{}'s \texttt{Utilities::MPI::some\_to\_some}. Then, each process $r$ merges all the points received from other processes into a single set $Y\src_r$ (\cref{alg:rbf-matrix-assembly}, line 13).

The evaluation points $\{\dstpoint{i}\}_{i = 1}^{N\dst}$ are treated in a similar way. The corresponding set of points owned by process $r$ is denoted by $X\dst_r$, and the set of points stored by process $r$ after inter-process communication is denoted by $Y\dst_r$. We remark that this procedure does not assume that the meshes $\mesh\src$ and $\mesh\dst$ are partitioned in a consistent way.

This step corresponds to lines 1--14 of \cref{alg:rbf-matrix-assembly}.

\begin{algorithm}[t]
  \caption{Parallel assembly of the matrices $\interpmatrix$ and $\evalmatrix$. The index $r$ denotes the rank of the current process in the MPI parallel communicator.}
  \label{alg:rbf-matrix-assembly}

  \begin{algorithmic}[1]
    \Require sets of locally owned points $X\src_r$ and $X\dst_r$
    \Require $h\src_{\max}$, maximum diameter of elements of $\mesh\src$
    \Ensure matrices $\interpmatrix$ and $\evalmatrix$
    \vspace{\baselineskip}

    \State construct a coarse representation $B_r$ of the local points $X\src_r$
    \State send $B_r$ to every parallel process different than $r$

    \State 
    \State $Y\src_{rr} \gets X\src_r$
    \State $Y\dst_{rr} \gets X\dst_r$
    \For{each parallel process $p$ different than $r$}
      \State $Y\src_{rp} \gets X\src_r \cap B_p$
      \State send the points $Y\src_{rp}$ to process $p$
      \State
      \State $Y\dst_{rp} \gets X\dst_r \cap B_p$
      \State send the points $Y\dst_{rp}$ to process $p$
    \EndFor
    \State $Y\src_r \gets \cup_{p = 1}^{N\text{proc}} Y\src_{pr}$
    \State $Y\dst_r \gets \cup_{p = 1}^{N\text{proc}} Y\dst_{pr}$

    \State 
    \For{every point $\srcpoint{j} \in X\src_r$}
      \State find $\mathbf z$, the $M$-th nearest neighbor (with respect to $\discgeodist$) to $\srcpoint{j}$ in $Y\src_r$
      \State $\rbfradius{j} \gets \max\left\{\alpha \discgeodist(\srcpoint{j}, \mathbf z), \rbfradiusmax{}\right\}$

      \State 
      \State compute the axis-aligned bounding box $C_j$ centered at $\srcpoint{j}$ and of dimension $2(\rbfradius{j} + 2h\src_{\max})$

      \State 
      \State $Z\src \gets Y\src_r \cap C_j$
      \For{every point $\srcpoint{i} \in Z\src$}
        \If{$\discgeodist(\srcpoint{j}, \srcpoint{i}) \leq \rbfradius{j}$}
          \State $(\interpmatrix)_{ij} \gets \rbf(d(\srcpoint{j}, \srcpoint{i}; \rbfradius{j}), \rbfradius{j})$
        \EndIf
      \EndFor

      \State 
      \State $Z\dst \gets Y\dst_r \cap C_j$
      \For{every point $\dstpoint{i} \in Z\dst$}
        \If{$\discgeodist(\srcpoint{j}, \dstpoint{i}) \leq \rbfradius{j}$}
          \State $(\evalmatrix)_{ij} \gets \rbf(d(\srcpoint{j}, \dstpoint{i}; \rbfradius{j}), \rbfradius{j})$
        \EndIf
      \EndFor
    \EndFor
  \end{algorithmic}
\end{algorithm}

\subsection{Finding the $M$-th nearest neighbor to a point}

In line 17 of \cref{alg:rbf-matrix-assembly}, we retrieve the $M$-th nearest neighbor to point $\srcpoint{j}$ within the set $Y\src_r$. This is achieved by iteratively building a binary max-heap \cite{williams1964algorithm} of $M$ elements. The point $\mathbf{z}$ is the root node of the heap. We remark that this task could be optimized by resorting to an M-tree data structure \cite{ciaccia1997indexing}.

Furthermore, the process of finding the $M$-th nearest neighbor is made faster by assuming that it is no further than a user-specified threshold $\rbfradiusmax{}$. While an optimal value for $\rbfradiusmax{}$ depends on $M$ and $\alpha$, we empirically selected it to be \num{10} times the average diameter $h\src_\text{avg}$ of elements of $\mesh\src$.

\subsection{Finding points sufficiently close to $\srcpoint{j}$}

Lines 20--22 of \cref{alg:rbf-matrix-assembly} have the purpose of quickly excluding the points of $Y\src_r$ that are sufficiently far from $\srcpoint{j}$. To do so, we observe that, with the same notation of \cref{sec:geodesic-distance-thresholding}, the following holds:
\begin{align}
  \|\mathbf{x} - \mathbf{y}\| &\leq \|\mathbf{x} - \mathbf{x}'\| + \|\mathbf{x}' - \mathbf{y}'\| + \|\mathbf{y} - \mathbf{y}'\| \\
  &\leq h\src_{\max} + \discgeodist(\mathbf x, \mathbf y) + h\src_{\max} \\
  &= \discgeodist(\mathbf x, \mathbf y) + 2h\src_{\max}\;,
\end{align}
where $h\src_{\max}$ is the maximum diameter of the elements of $\mesh\src$.
Now, let $C_j$ be an axis-aligned bounding box centered at $\srcpoint{j}$ of dimension $2(r_j + 2h\src_{\max})$ along all coordinate directions (line 20). If a point $\srcpoint{i}$ lies outside the box, there holds
\begin{gather}
  \|\srcpoint{j} - \srcpoint{i}\| \geq r_j + 2h\src_{\max}\;, \\
  \discgeodist(\srcpoint{j}, \srcpoint{i}) \geq \|\srcpoint{j} - \srcpoint{i}\| - 2h\src_{\max} \geq r_j\;,
\end{gather}
so that the associated entry $(\interpmatrix)_{ij}$ in the interpolation matrix is zero. Consequently, we can discard all points that lie outside of the bounding box, without explicitly computing their discrete geodesic distance from $\srcpoint{j}$.

This is especially convenient because the intersection operation of line 22 can be performed very efficiently by storing the points $Y\src_r$ into an R-tree structure \cite{manolopoulos2006rtrees, guttman1984rtrees}. We rely on the R-tree implementation offered by the Boost library \cite{boost}.

\subsection{Evaluating the discrete geodesic distance}
\label{sec:algorithm-geodist}

We use Dijkstra's shortest path algorithm \cite{dijkstra1959note} to evaluate $\discgeodist$, as described in \cref{alg:dijkstra}. As usually done, the set of visited points $A$ is stored in a binary min-heap data structure \cite{williams1964algorithm}, providing constant-time retrieval of the element with minimum tentative distance (\cref{alg:dijkstra}, line 9).

Since we only use $\discgeodist$ to threshold the Euclidean distance in \eqref{eq:distance-cutoff}, we stop the algorithm if it reaches nodes that are further than the threshold $R$ from the starting point (\cref{alg:dijkstra}, line 21). Similarly, if the two endpoints are further than $R$ in Euclidean distance, we do not compute their discrete geodesic distance (\cref{alg:dijkstra}, lines 1--3). In both cases, we return an arbitrary value greater than or equal to $R$.

The loop of \cref{alg:rbf-matrix-assembly}, lines 16--35, requires to evaluate $\discgeodist(\srcpoint{j}, \mathbf y)$ multiple times with $\srcpoint{j}$ as first argument and varying the second argument $\mathbf y$. To avoid unnecessary computations, the sets $A$ and $B$ and the tentative distances $s_\mathbf{z}$ are reused between calls. If the point $\mathbf y$ was already added to the set $B$ during a previous evaluation of the distance, we reuse the already computed distance $s_\mathbf{y}$.

When running in parallel, the mesh $\mesh$ is distributed over multiple processes, and each of the processes stores only a subset of the mesh vertices. However, each process needs to evaluate the discrete geodesic distance between arbitrary mesh vertices, regardless of whether they are owned by that process or not. A naive solution would be to send all the vertices and their adjacency information to every process. This however does not scale efficiently over large meshes or large number of processes.

To avoid this issue, we observe that, when evaluating $\discgeodist$, each process $r$ always passes a locally owned point $\srcpoint{j} \in X\src_r$ as first argument (see \cref{alg:rbf-matrix-assembly}, lines 18, 24 and 31). Furthermore, we are only interested in destination points no further than the maximum \ac{RBF} support radius $\rbfradiusmax{}$. Therefore, it is sufficient that each process stores the points that are within $\rbfradiusmax{}$ of any locally owned point, together with the associated adjacency information.

\begin{algorithm}[t]
  \caption{Dijkstra's shortest path algorithm, adapted from \cite{dijkstra1959note}.}
  \label{alg:dijkstra}

  \begin{algorithmic}[1]
    \Require reference mesh $\mesh$
    \Require start and end mesh vertices $\mathbf x$, $\mathbf y$
    \Require threshold $R$
    \Ensure returns $\discgeodist(\mathbf x, \mathbf y)$
    \vspace{\baselineskip}
    \If{$\|\mathbf x - \mathbf y\| > R$}
      \State \textbf{return} $R$
    \EndIf
    \State
    \State add $\mathbf x$ to the set $A$ of visited points
    \State set the tentative distance for $\mathbf x$: $s_{\mathbf x} \gets 0$
    \State
    \While{$A$ is not empty}
      \State find the point $\mathbf z$ in $A$ for which the tentative distance $s_{\mathbf z}$ is minimum
      \State remove $\mathbf z$ from $A$ and add it to the set of closed points $B$
      \State
      \For{each neighbor $\mathbf w$ of $\mathbf z$ that is not in $B$}
        \If{$\mathbf w \in A$}
          \State $s_{\mathbf w} \gets \min\{s_{\mathbf z} + \|\mathbf w - \mathbf z\|, s_{\mathbf w}\}$
        \Else
          \State add $\mathbf w$ to $A$
          \State $s_{\mathbf w} \gets s_{\mathbf z} + \|\mathbf w - \mathbf z\|$
        \EndIf
      \EndFor
      \State
      \If{$\mathbf z = \mathbf y$ \textbf{or} $s_{\mathbf z} > R$}
        \State \textbf{return} $s_{\mathbf z}$
      \EndIf
    \EndWhile
  \end{algorithmic}
\end{algorithm}

\subsection{Evaluating the interpolant}

After constructing the matrices $\interpmatrix$ and $\evalmatrix$, given an input data vector $\mathbf{f}$, the interpolation coefficients are obtained by solving the linear system \eqref{eq:rbf-linear-system} with the preconditioned GMRES method \cite{saad2003iterative}. For the numerical experiments of this paper, we consider an \ac{ILU} factorization preconditioner with \num{2} levels of fill-in and \num{1} level of inter-process overlap, as implemented by Trilinos IFPACK \cite{sala2005robust} and wrapped by \dealii{}. Alternatively, one can consider an ad-hoc preconditioner based on approximate cardinal basis functions \cite{bucelli2023preserving,brown2005approximate,beatson1999fast}, although in our numerical experiments it performed less effectively than the \ac{ILU} preconditioner.

\begin{figure}
  \centering

  \begin{subfigure}{0.24\textwidth}
    \includegraphics[width=\textwidth]{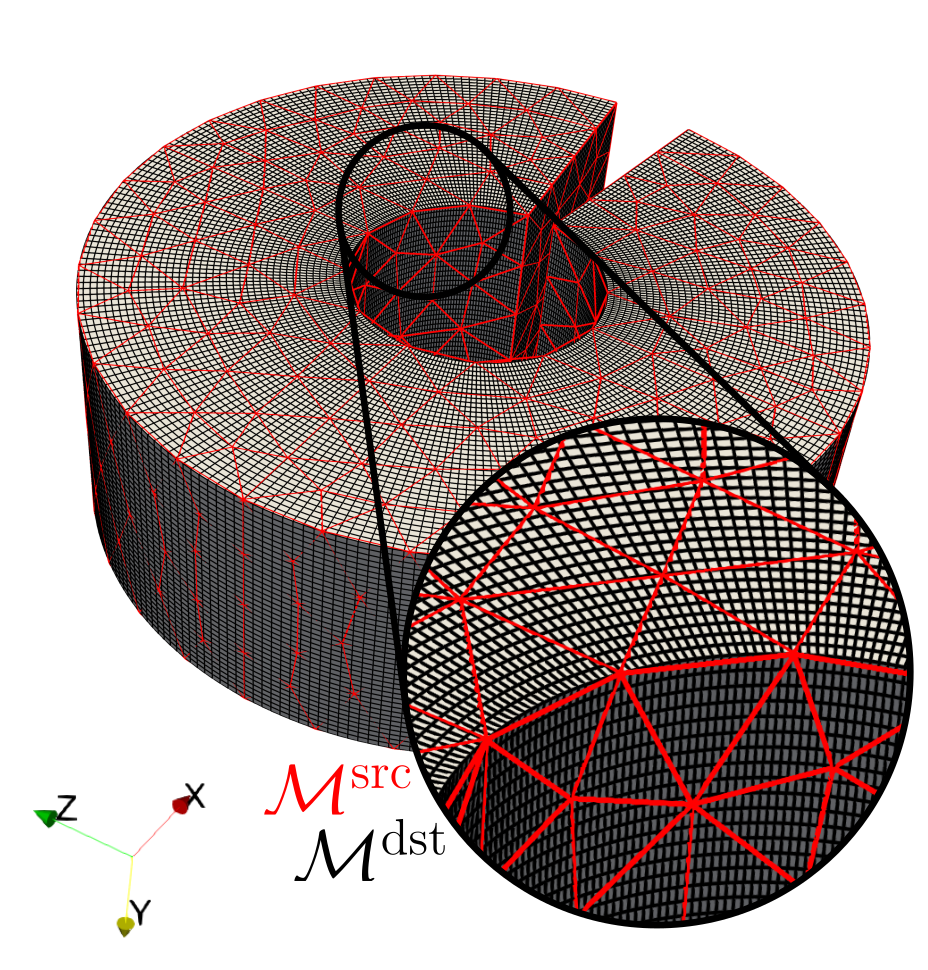}
    \caption{}
    \label{fig:ring-domain}
  \end{subfigure}
  \begin{subfigure}{0.24\textwidth}
    \includegraphics[width=\textwidth]{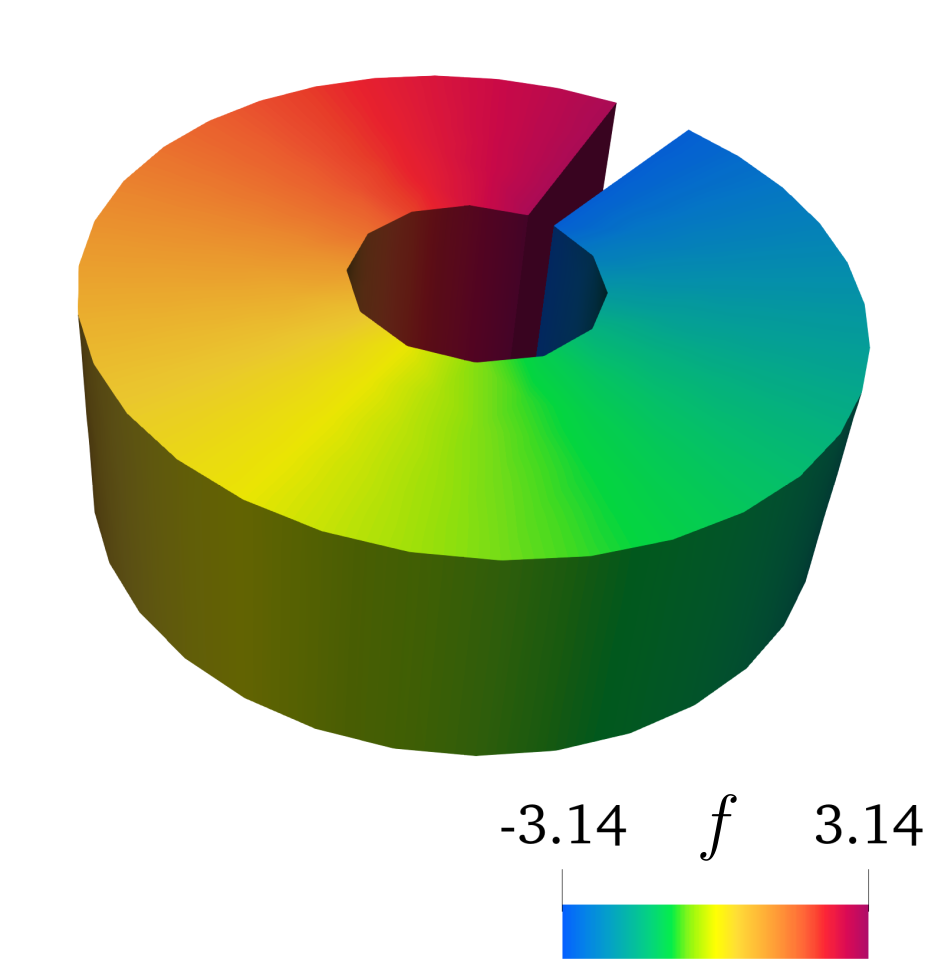}
    \caption{source data}
    \label{fig:ring-source-data}
  \end{subfigure}
  \begin{subfigure}{0.24\textwidth}
    \includegraphics[width=\textwidth]{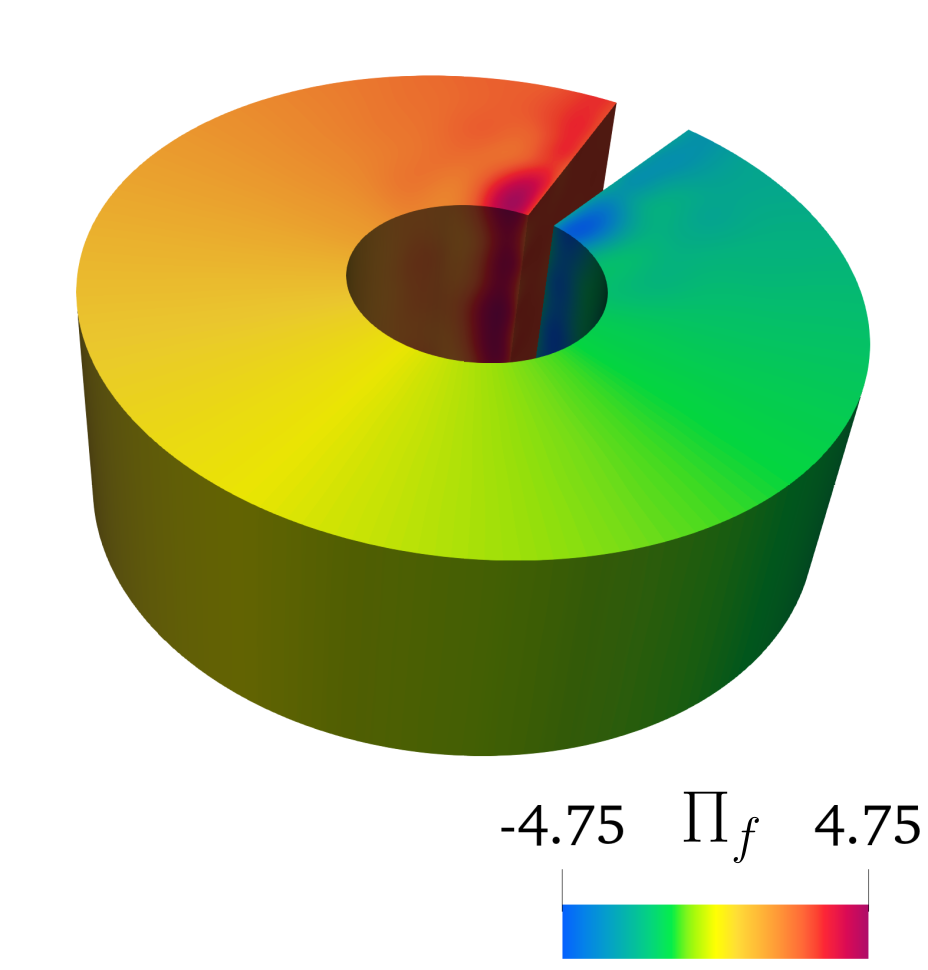}
    \caption{\ac{RLRBF}}
    \label{fig:ring-destination-nothreshold}
  \end{subfigure}
  \begin{subfigure}{0.24\textwidth}
    \includegraphics[width=\textwidth]{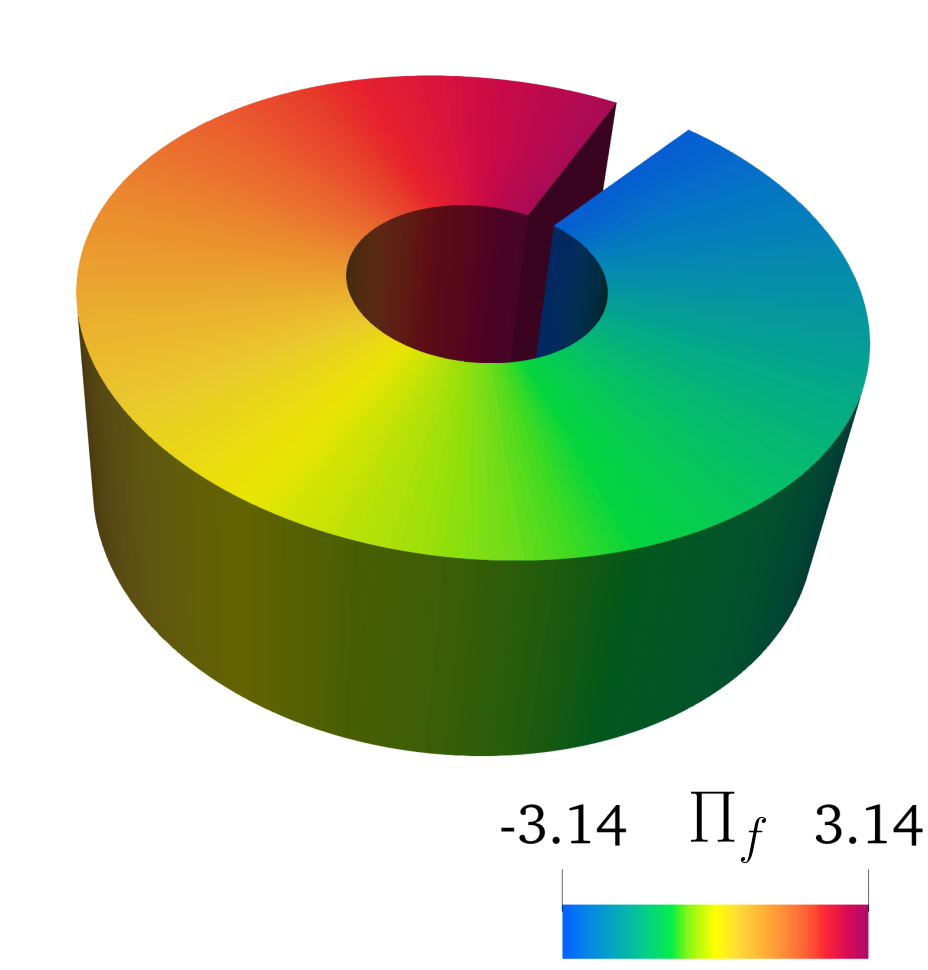}
    \caption{\rlrbfg{}}
    \label{fig:ring-destination-threshold}
  \end{subfigure}

  \caption{Idealized benchmark (\cref{sec:ring}). (a) Interpolation domain $\Omega$, overlaid with the meshes $\mesh\src = \mesh^\text{tet}_{\text{ring,1}}$ (red) and $\mesh\dst = \mesh^\text{hex}_{\text{ring,1}}$ (black). (b) Source data $f(x, y, z) = \atantwo(z, -x)$. (c) Interpolant obtained without geodesic thresholding. Notice the different colorbar scale and the spurious oscillations close to the slit. (d) Interpolant obtained from mesh $\mesh^\text{tet}_{\text{ring,1}}$ with geodesic thresholding.}
\end{figure}

\section{Numerical experiments}
\label{sec:numerical-experiments}

In the following sections, we examine the properties of the proposed interpolation method by means of numerical experiments, and compare it with \ac{RLRBF} interpolation. Unless otherwise specified, we set $\beta = \num{0.5}$. The role of the high-curvature threshold is further discussed in \cref{sec:high-curvature}.

\subsection{Idealized benchmark}
\label{sec:ring}

We consider an idealized benchmark on the domain depicted in \cref{fig:ring-domain}. The domain represents a toroid with square section, from which a thin slice has been removed, opening a slit. We interpolate the following function:
\begin{equation}
  f(x, y, z) = \atantwo(z, -x) = \begin{dcases}
    \arctan(-\frac{z}{x}) & \text{if } x < 0\;, \\
    \arctan(-\frac{z}{x}) + \pi & \text{if } x > 0 \text{ and } z \geq 0\;, \\
    \arctan(-\frac{z}{x}) - \pi & \text{if } x > 0 \text{ and } z < 0\;, \\
    \frac{\pi}{2} & \text{if } x = 0 \text{ and } z > 0\;, \\
    -\frac{\pi}{2} & \text{if } x = 0 \text{ and } z < 0\;, \\
    \text{undefined} & \text{if } x = z = 0\;,
  \end{dcases}
\end{equation}
which is discontinuous across the slit, as shown in \cref{fig:ring-source-data}. We interpolate the function $f$ as described in previous sections, both with and without geodesic distance thresholding. \Cref{tab:mesh} lists the meshes considered throughout this section.

\Cref{fig:ring-destination-nothreshold} shows an example of \ac{RLRBF} interpolant, using the mesh $\mesh^\text{tet}_{\text{ring,1}}$ and setting $M = \num{6}$, $\alpha = \num{2}$. Close to the slit, the interpolant features large oscillations, significantly departing from the original function $f$. In intuitive terms, this is due to $f$ having very different values on the two sides of the slit. Without geodesic thresholding, the interpolant does not take the slit into account, and the oscillations result from trying to fit the values on both sides. On the contrary, constructing the interpolant with geodesic thresholding, with the same mesh and parameters, leads to a much more accurate result, as shown in \cref{fig:ring-destination-threshold}.

\begin{table}
  \centering

  \begin{tabular}{c c c S S[table-format=2.2,table-auto-round] S[table-format=2.2,table-auto-round] S[table-format=2.2,table-auto-round]}
    \toprule
    \textbf{Domain} & \textbf{Mesh} & \textbf{Element type} & \textbf{\# vertices} & $h_{\min}$ [\si{\milli\metre}] & $h_\text{avg}$ [\si{\milli\metre}] & $h_{\max}$ [\si{\milli\metre}] \\
    \midrule
    ring & $\mesh^\text{tet}_{\text{ring,1}}$ & tet & 608 & 25 & 40.0225 & 55.8842 \\
    ring & $\mesh^\text{tet}_{\text{ring,2}}$ & tet & 1166 & 19.4756 & 31.2794 & 41.7913 \\
    ring & $\mesh^\text{tet}_{\text{ring,3}}$ & tet & 2453 & 13.412 & 23.3047 & 31.1178 \\
    ring & $\mesh^\text{tet}_{\text{ring,4}}$ & tet & 6403 & 9.19972 & 16.277 & 22.6551 \\
    ring & $\mesh^\text{tet}_{\text{ring,5}}$ & tet & 14957 & 7.14286 & 11.8231 & 15.3818 \\
    ring & $\mesh^\text{tet}_{\text{ring,6}}$ & tet & 42327 & 4.86997 & 8.14003 & 11.0399 \\
    ring & $\mesh^\text{tet}_{\text{ring,7}}$ & tet & 304313 & 2.47238 & 4.08136 & 5.64289 \\
    ring & $\mesh^\text{tet}_{\text{ring,8}}$ & tet & 4411070 & 0.994159 & 1.6398 & 2.3575 \\
    ring & $\mesh^\text{hex}_{\text{ring}}$ & hex & 418569 & 3.73385 & 4.33481 & 5.16729 \\
    \midrule
    heart & $\mesh^\text{fine}_\text{heart}$ & tet & 1260388 & 0.36081 & 0.865568 & 3.11138 \\
    heart & $\mesh^\text{coarse}_\text{heart}$ & tet & 67553 & 0.821286 & 2.50623 & 4.80932 \\
    \bottomrule
  \end{tabular}

  \caption{Meshes used in the numerical experiments. For each mesh, we report the type of elements (hexahedra or tetrahedra), the number of nodes, and the minimum, average and maximum element diameter $h$.}
  \label{tab:mesh}
\end{table}

\subsubsection{Convergence test}
\label{sec:convergence}

We assess the convergence of the proposed method by constructing the interpolant of $f$ on meshes $\mesh^\text{tet}_{\text{ring},i}$, for $i = 1, 2, \dots, 7$, and evaluating it on $\mesh^\text{hex}_{\text{ring},1}$ (see \cref{tab:mesh}), both with and without geodesic thresholding. We measure the error between the interpolant and $f$ with the following approximation of the relative $L^\infty(\Omega)$ norm:
\begin{equation}
  \errinfty = \frac{\max_{i = 1, \dots, N\dst}\left|\mathbf f\dst_i - f(\dstpoint{i})\right|}{\max_{i = 1, \dots, N\dst}\left|f(\dstpoint{i})\right|}\;.
  \label{eq:error}
\end{equation}

\begin{figure}
  \centering
  \includegraphics{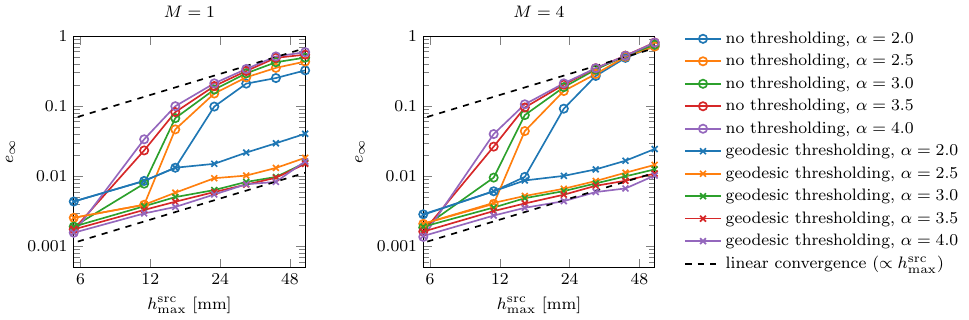}

  \caption{Idealized benchmark, convergence test (\cref{sec:convergence}). Interpolation error, defined as in \eqref{eq:error}, against the maximum element diameter $h\src_{\max}$ of the source mesh $\mesh\src$. The interpolant was computed setting $M = 1$ (left) and $M = 4$ (right), for different values of $\alpha$. Different choices of $\alpha$ correspond to different colors, while circles and crosses correspond to interpolants without and with geodesic thresholding, respectively.}
  \label{fig:convergence}
\end{figure}

\Cref{fig:convergence} reports the error obtained for $M = 1$ and $M = 4$ with different choices of $\alpha$, both with and without geodesic thresholding. For relatively large values of $h\src_{\max}$, \rlrbfg{} interpolation reduces the interpolation error approximately by one order of magnitude with respect to \ac{RLRBF}. Indeed, the latter is largely affected by the oscillations appearing near the slit.

The coefficient $\alpha$ rescales the \ac{RBF} support radius. For large values of $h\src_{\max}$, increasing $\alpha$ without geodesic thresholding allows points on one side of the slit to interact more with points on the other, leading to larger errors. Conversely, increasing $\alpha$ with geodesic thresholding reduces the error, since a larger \ac{RBF} support radius removes small spurious oscillations in the interpolant, improving its accuracy.

For a sufficiently small $h\src_{\max}$, \ac{RLRBF} and \rlrbfg{} lead to approximately equal interpolation errors. This happens because the interpolation radii $\rbfradius{j}$ are progressively reduced as the mesh is refined, so that for a sufficiently fine mesh the interpolation radius is small enough that points on one side of the slit do not interact with points on the other side. In that case, geodesic thresholding has no effect in practice, and the two methods yield almost identical results.

Finally, our numerical results indicate that the \rlrbfg{} interpolant converges to $f$ with approximate order \num{1} with respect to $h\src_{\max}$. This is in agreement with the analysis presented in \cite{demarchi2020convergence} for \ac{RLRBF} interpolation.

\subsubsection{Computational efficiency and scalability}
\label{sec:scalability}

In the same setting as previous section, we run a strong scalability test by computing the interpolant from mesh $\mesh^\text{tet}_{\text{ring},8}$ onto mesh $\mesh^\text{hex}_{\text{ring}}$ varying the number of parallel processes. We measure the wall time needed for the different stages of the construction of the interpolant and its evaluation. The scalability test was run on the GALILEO100 supercomputer\footnote{Technical specifications: \url{https://wiki.u-gov.it/confluence/display/SCAIUS/UG3.3\%3A+GALILEO100+UserGuide}.} from the CINECA supercomputing center (Italy).

\begin{figure}
  \centering
  \includegraphics{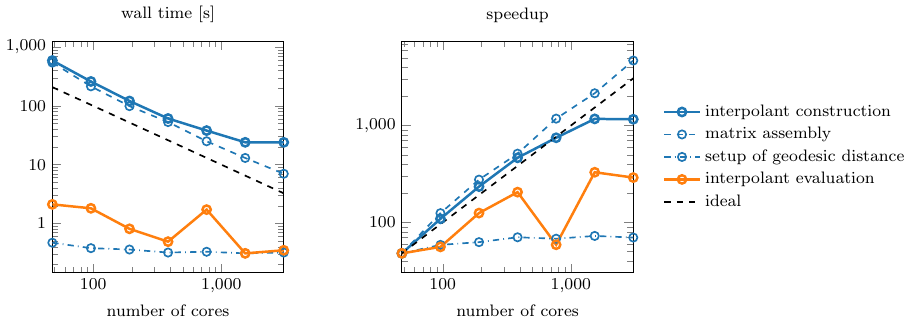}

  \caption{Idealized benchmark, scalability test (\cref{sec:scalability}). Total wall-time (left) and parallel speedup for different steps of the construction and evaluation of the \rlrbfg{} interpolant.}
  \label{fig:scalability}
\end{figure}

The results are shown in \cref{fig:scalability}, where we report the computational time associated with the wall time taken by the interpolant construction (which includes the initialization of the adjacency graph for evaluating the discrete geodesic distance, the construction of the matrices $\interpmatrix$ and $\evalmatrix$, the initialization of the preconditioner for the system \eqref{eq:rbf-linear-system}, and the evaluation of $\interp{g}$ for rescaling) and the time needed for the evaluation of the interpolant (i.e. for the solution of the linear system \eqref{eq:rbf-linear-system}). We also report separately the time to assemble $\interpmatrix$, $\evalmatrix$ and the preconditioner (``matrix assembly'' in the plot) and the time required to initialize the adjacency graph for geodesic distance (``setup of geodesic distance'' in the plot).

The time needed to assemble matrices scales ideally up to \num{3072} cores. The scalability of the overall interpolant construction is ideal up to around \num{768} cores (corresponding to approximately \num{5500} interpolation points per core), then starts to deteriorate due to parallel communication and load balancing issues.

The solution of the interpolation system \eqref{eq:rbf-linear-system} exhibits relatively poor scalability properties, mostly due to the scalability of the preconditioner. Furthermore, the time needed for the initialization of the discrete geodesic distance data structures is essentially independent of the number of cores. Nonetheless, both these steps have a negligible impact on the overall computational cost.

\subsection{Cardiac electromechanics test}
\label{sec:heart}

With the aim of providing an example of application to a challenging multiphysics problem of practical interest, we consider a numerical model of whole-heart cardiac electromechanics described in \cite{fedele2023comprehensive}. The term cardiac electromechanics refers to the interplay of electrical excitation (known as cardiac electrophysiology, \cite{collifranzone2014mathematical, piersanti2021modeling, sundnes2007computing, sung2020personalized, gillette2021framework}), active contraction and passive mechanical response of the cardiac muscle during the heartbeat \cite{quarteroni2019mathematical,augustin2021computationally,gurev2011models,stella2022fast,regazzoni2022cardiac,piersanti20223d,fedele2023comprehensive,strocchi2020simulating,gerach2024differential,quarteroni2022modeling,quarteroni2023mathematical}. Such models have the potential of supporting the understanding and clinical decision making associated with cardiovascular pathologies \cite{salvador2022role, peirlinck2021precision, salvador2021electromechanical}.

Cardiac electromechanics is intrinsically multiphysics and multiscale. Indeed, electrophysiology is characterized by much smaller spatial and temporal scales than muscular mechanics (both active and passive). For this reason, the electrophysiology model is typically solved on a finer discretization than mechanics, both in time and in space, by combining temporal staggering and intergrid transfer operators \cite{regazzoni2022cardiac, piersanti20223d, fedele2023comprehensive, salvador2020intergrid, bucelli2023preserving, gerach2021electro, viola2023gpu}.

The heart is composed of four different interconnected chambers, and it is connected to large blood vessels (including the aorta, the pulmonary artery), resulting in a topologically complex domain (see e.g. \cref{fig:complex-domains}). It is therefore crucial to use intergrid transfer operators that account for this complexity. To this end, we apply the newly introduced \rlrbfg{} method in the electromechanical framework.

\begin{figure}
  \centering

  \includegraphics[width=0.8\textwidth]{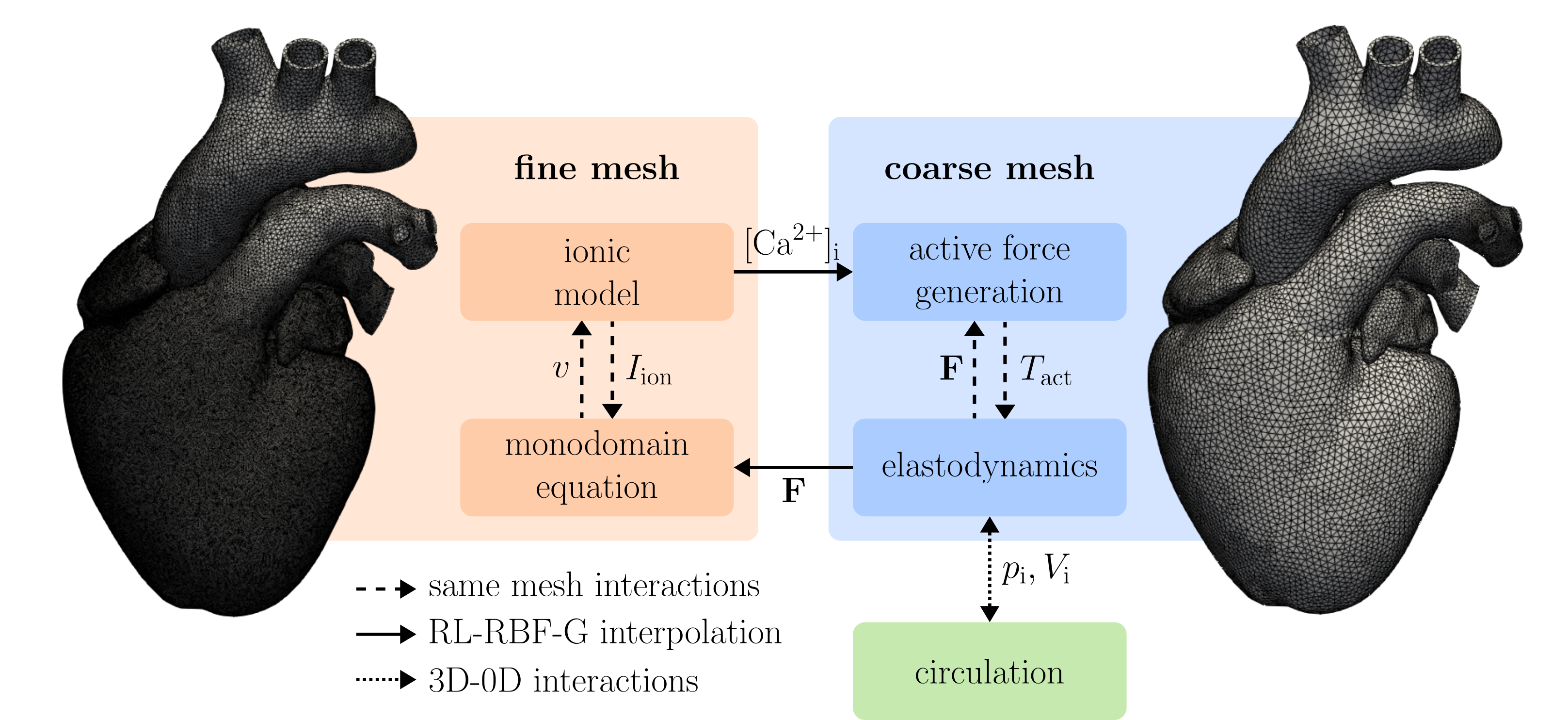}
  \caption{Schematic representation of the electromechanics model of \cref{sec:heart}, reporting the multiphysics interactions.}
  \label{fig:em-model}
\end{figure}

The whole-heart electromechanical model features the following components (see also \cref{fig:em-model}) \cite{fedele2023comprehensive}:
\begin{itemize}
  \item the monodomain model of cardiac electrophysiology, for the spatial propagation of electrical excitation \cite{collifranzone2014mathematical}. The model includes geometry-dependent \ac{MEF} effects: denoting by $\displacement$ the muscular displacement, the monodomain model accounts for the effect of the deformation onto the electrical conductivity, through $\Ftensor = \Itensor + \grad\displacement$ and $\jacobian = \det\Ftensor$ \cite{salvador2022role, bucelli2023preserving, gerach2024differential};
  \item the \ac{TTP} and \ac{CRN} ionic models, describing the electrical activity at the cellular scale \cite{collifranzone2014mathematical, ten2006alternans,courtemanche1998ionic}. The ionic models regulate in particular the evolution of the intracellular calcium concentration $\calcium$;
  \item the RDQ20-MF model of muscular contraction, that computes an active stress tensor $\actstress$ in response to variations in $\calcium$ \cite{regazzoni2020biophysically};
  \item the elastodynamics equation for muscular deformation, including both passive and active stress contributions \cite{regazzoni2022cardiac};
  \item a lumped-parameter model of the circulatory system, describing pressures and flowrates through several circulation compartments by means of a system of \acp{ODE} \cite{hirschvogel2017monolithic, blanco20103d, fedele2023comprehensive}.
\end{itemize}
The model equations are reported in full in \ref{app:em-model}. We refer to \cite{fedele2023comprehensive} for further details on the electromechanical model and strategies for its numerical discretization, and to \cite{zingaro2023electromechanics} for a discussion on its calibration.

Similarly to \cite{regazzoni2022cardiac, piersanti20223d, bucelli2023preserving, salvador2020intergrid}, we solve the ionic and monodomain models on a fine mesh $\mesh^\text{fine}_\text{heart}$, and the force generation and mechanics models on a coarser mesh $\mesh^\text{coarse}_\text{heart}$ (see \cref{tab:mesh}). We interpolate the calcium concentration $\calcium$, computed on the fine mesh by the ionic model, onto the coarse mesh, where it is needed by the force generation model. Conversely, we interpolate the deformation gradient $\Ftensor$ computed by solving the elastodynamics equation from the coarse to the fine mesh, using the method presented in \cite{bucelli2023preserving}, wherein we replace the \ac{RLRBF} interpolation with the newly introduced \rlrbfg{} method. For the sake of post-processing, we also interpolate the displacement field $\displacement$ from the coarse to the fine mesh.

We run the simulation up to a final time of $T = \SI{4.8}{\second}$, corresponding to a total of six heartbeats, with each heartbeat lasting $t_\text{HB} = \SI{0.8}{\second}$. We report results for the last hearbeat only (from $t = \SI{4.0}{\second}$ to $t = \SI{4.8}{\second}$). We consider the segregated-staggered scheme previously discussed in \cite{regazzoni2022cardiac, piersanti20223d, fedele2023comprehensive}. We set a time discretization step of $\Delta t_\text{EP} = \SI{5e-5}{\second}$ for the monodomain and ionic models, while force generation, muscular mechanics and circulation are solved with a time discretization step of $\Delta t_\text{M} = 20\Delta t_\text{EP} = \SI{1e-3}{\second}$. The simulation is run using \num{92} cores from the iHEART computing node\footnote{The node runs eight processors Intel Xeon Platinum 8160 @\SI{2.1}{\giga\hertz}.} available at MOX, Dipartimento di Matematica, Politecnico di Milano.

\begin{figure}
  \centering

  \includegraphics[width=0.85\textwidth]{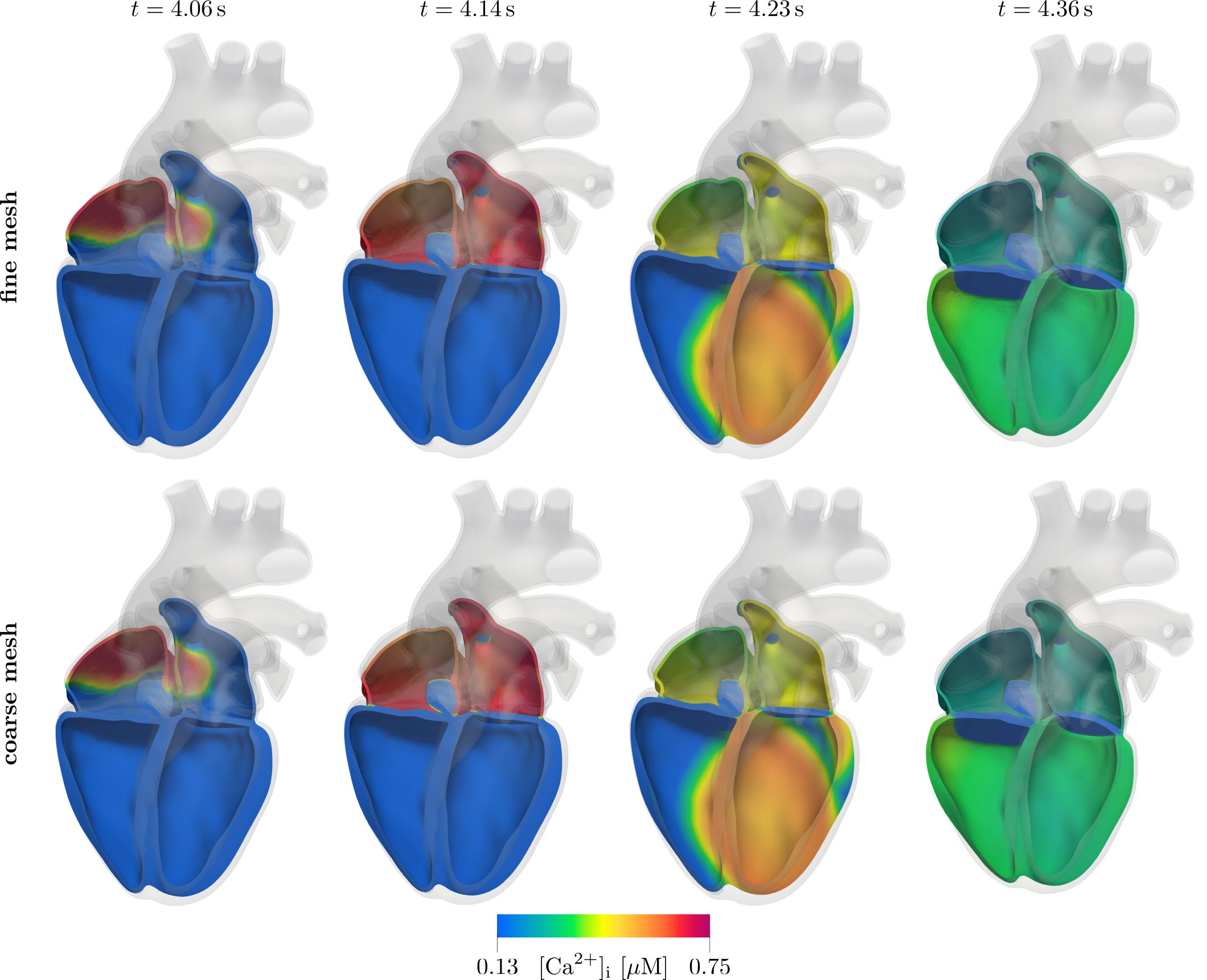}

  \caption{Whole-heart electromechanics (\cref{sec:heart}). Intracellular calcium concentration $\calcium$ over time on the electrophysiology (top) and mechanics (bottom) meshes. The domain is warped according to the displacement $\displacement$.}
  \label{fig:em-calcium}
\end{figure}

\begin{figure}
  \centering

  \includegraphics[width=0.85\textwidth]{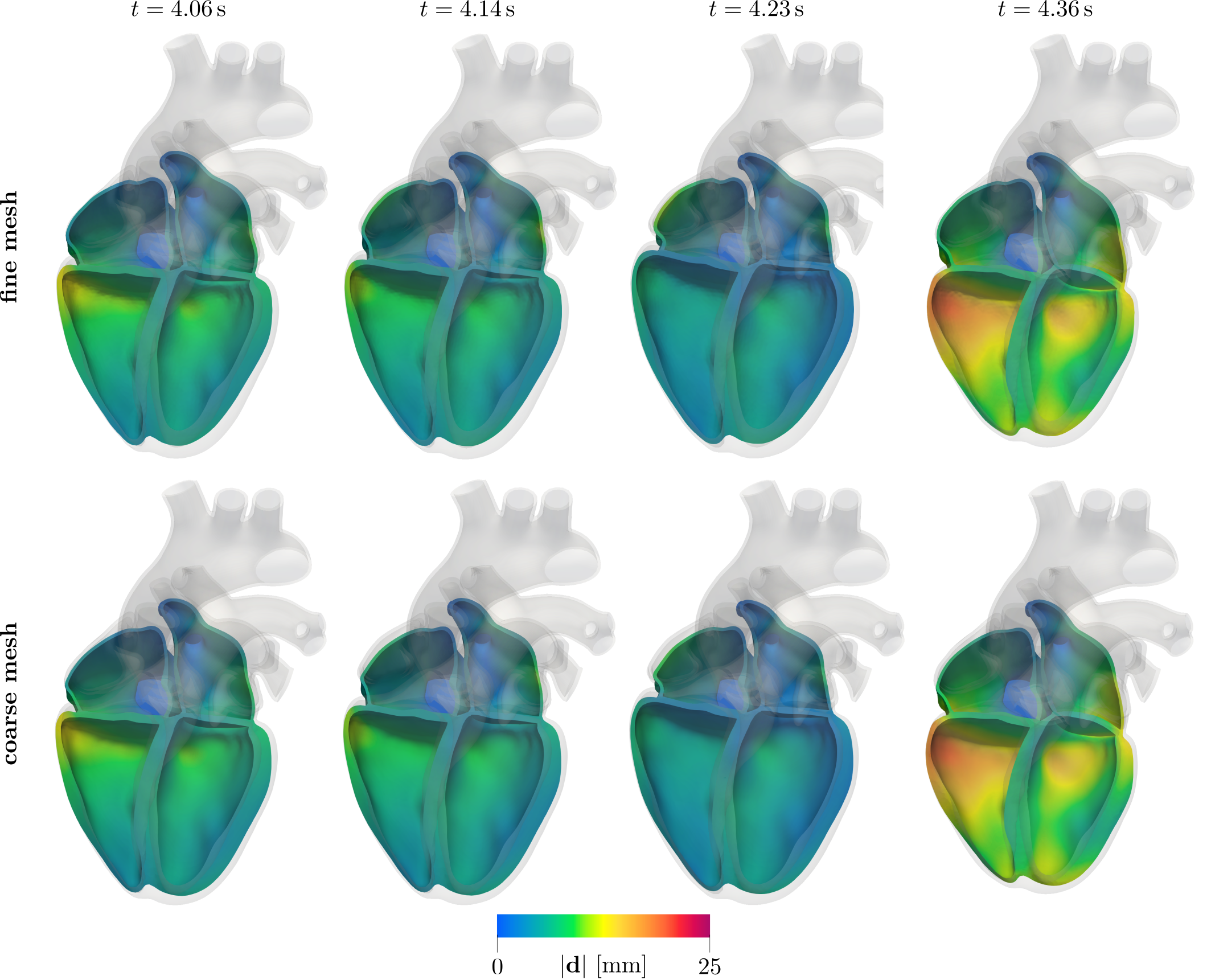}

  \caption{Whole-heart electromechanics (\cref{sec:heart}). Muscular displacement $\displacement$ over time on the electrophysiology (top) and mechanics (bottom) meshes. We remark that the displacement is interpolated for post-processing, while the monodomain equation relies on interpolating the deformation gradient $\Ftensor$. The domain is warped according to $\displacement$.}
  \label{fig:em-displacement}
\end{figure}

\begin{figure}
  \centering

  \includegraphics{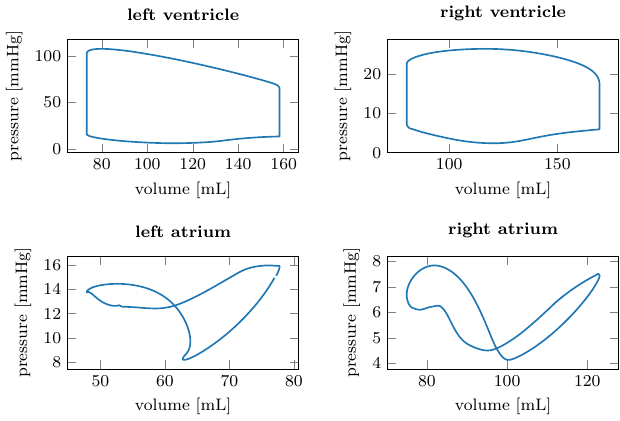}

  \caption{Whole-heart electromechanics (\cref{sec:heart}). Pressure-volume loops of the four cardiac chambers.}
  \label{fig:em-pvloop}
\end{figure}

The two meshes have a significantly different size, with the fine mesh having almost \num{20} times as many nodes as the coarse one. Nonetheless, our numerical results indicate that the \rlrbfg{} method allows to reliably transfer information between the two meshes. No spurious oscillations are observed in the interpolated data, and the overall results are in agreement with those of \cite{fedele2023comprehensive}, as seen in \cref{fig:em-calcium,fig:em-displacement,fig:em-pvloop}, reporting the evolution of $\calcium$, $\displacement$ and the pressure-volume loops of the four cardiac chambers.

\subsubsection{Computational cost of interpolation}

\Cref{tab:heart-cost} reports a breakdown of the computational cost associated with the simulation, using either \rlrbfg{} or \ac{RLRBF} interpolations. The assembly of interpolation matrices takes up approximately \SI{17}{\percent} of the total cost of the initialization (which is mostly spent in computing the initial conditions for the mechanics problem). Interpolating the deformation gradient tensor $\Ftensor$ from the coarse to the fine mesh requires approximately \SI{12}{\percent} of the average time needed to solve one electromechanics time step (comprising \num{20} time steps of electrophysiology, one time step of force generation and one of mechanics, one fine-to-coarse interpolation and one coarse-to-fine interpolation). We remark that, as described in \cite{bucelli2023preserving}, this entails the interpolation of \num{11} scalar fields that define the tensor $\Ftensor$. Finally, the interpolation of $\calcium$ from the coarse to the fine mesh has a very small impact on the cost of one electromechanics time step.

The cost of initializing \ac{RLRBF} interpolation is significantly lower than the one of \rlrbfg{}. This is largely due to the need to compute geodesic distances in constructing the interpolation matrices with \rlrbfg{}. Furthermore, evaluating the interpolant of $\Ftensor$ with \rlrbfg{} is more costly than with \ac{RLRBF}, due to a large number of linear solver iterations in solving \eqref{eq:rbf-linear-system}. Nonetheless, the large oscillations appearing in the \ac{RLRBF} interpolant prevent its use for this application.

Although not negligible, the overall cost of constructing and evaluating the \rlrbfg{} interpolants does not dominate the wall time of the simulations. We believe this overhead is justified by the discretization flexibility and the reliability that the proposed interpolation method offers. Indeed, the \rlrbfg{} interpolation allows to model cardiac electrophysiology to high accuracy without impacting the computational cost associated with the mechanics submodel (contrarily to the strategy of \cite{fedele2023comprehensive}). Therefore, it can support simulations where the electrophysiological accuracy is critical, e.g. including spatial heterogeneities and electrical dysfunctions \cite{salvador2021electromechanical, pagani2021computational, heijman2021computational}.

\begin{table}
  \centering

  \begin{tabular}{l
                  S[table-format=4.2,table-alignment-mode=none]
                  S[table-format=3.1,table-auto-round]
                  S[table-format=4.2,table-alignment-mode=none]
                  S[table-format=3.1,table-auto-round]}
    \toprule
    & \multicolumn{2}{c}{\textbf{\rlrbfg{}}} & \multicolumn{2}{c}{\textbf{\ac{RLRBF}}} \\
    & \textbf{wall time } [\si{\second}] & \textbf{relative } [\si{\percent}] & \textbf{wall time } [\si{\second}] & \textbf{relative } [\si{\percent}] \\
    \midrule
    \textbf{simulation setup}           & 2940 & 100.00 & 2300 & 100.00 \\
    \qquad mechanics                    & 2129 &  72.41 & 1994 &  86.73 \\
    \qquad fine-to-coarse interpolant   &  266 &   9.05 &   23 &   1.01 \\
    \qquad coarse-to-fine interpolant   &  226 &   7.69 &   50 &   2.18 \\
    \qquad electrophysiology            &  108 &   3.67 &   94 &   4.08 \\
    \qquad fibers                       &   17 &   0.58 &   16 &   0.68 \\
    \qquad active stress                &  0.1 & {<0.01}&  0.1 & {<0.01}\\
    \midrule
    \textbf{electromechanics time step} & 18.00 & 100.00 & 16.80 & 100.00 \\
    \qquad electrophysiology            &  6.24 &  34.66 &  6.22 &  37.03 \\
    \qquad mechanics                    &  9.21 &  51.17 &  9.23 &  54.96 \\
    \qquad coarse-to-fine interpolation &  2.18 &  12.11 &  0.98 &   5.86 \\
    \qquad fine-to-coarse interpolation &  0.37 &   2.06 &  0.36 &   2.16 \\
    \bottomrule
  \end{tabular}

  \caption{Breakdown of the computational cost for the whole-heart electromechanics simulation (\cref{sec:heart}), using \rlrbfg{} or \ac{RLRBF} interpolation. We report the wall time (WT) associated with the simulation setup and with solving one time step of electromechanics. We also report the time spent in substeps, with the ``relative'' column reporting the percentage of the substep time relative to the parent step.}
  \label{tab:heart-cost}
\end{table}

\subsubsection{Role of the high-curvature threshold}
\label{sec:high-curvature}

To better understand the role of the high-curvature detection coefficient $\beta$, we perform a test using the same setting of previous section, except that we set $\beta = +\infty$ in \eqref{eq:distance-cutoff}. This is equivalent to disabling the detection of high-curvature regions altogether, i.e. redefining
\begin{equation}
  d(\mathbf x, \mathbf y; R) = \begin{cases}
    +\infty & \text{if } \discgeodist(\mathbf x, \mathbf y) > R\;, \\
    \|\mathbf x - \mathbf y\| & \text{otherwise.}
  \end{cases}
\end{equation}

\begin{figure}
  \centering
  \includegraphics[width=0.85\textwidth]{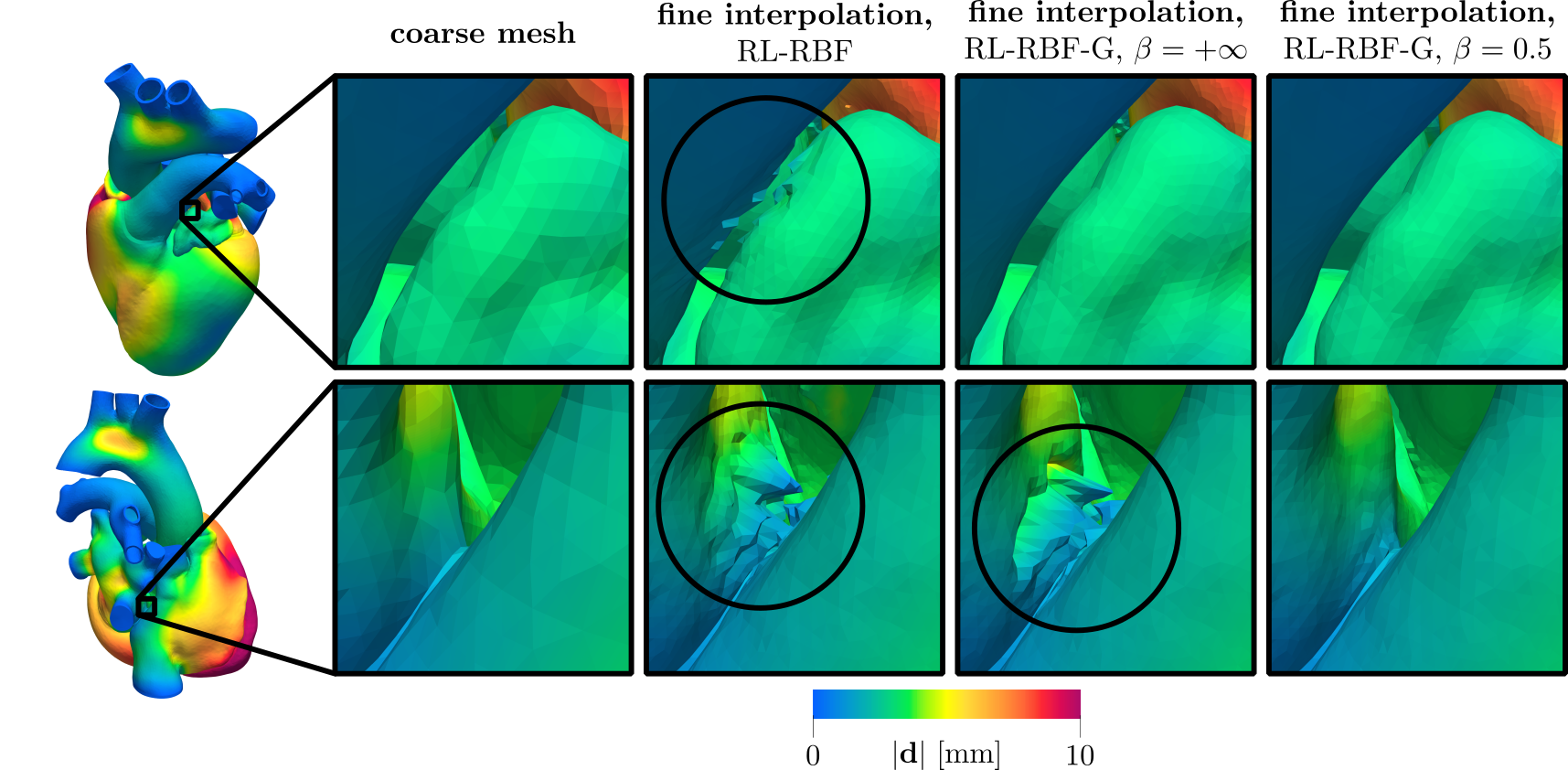}

  \caption{Whole-heart electromechanics (\cref{sec:heart}). Effect of \rlrbfg{} interpolation and of the high-curvature detection in two regions close to the left atrial appendage and the fossa ovalis. Circles indicate spurious oscillations.}
  \label{fig:high-curvature}
\end{figure}

\Cref{fig:high-curvature} highlights two critical regions for the interpolation procedure, close to the left atrial appendage and to the fossa ovalis i.e. a region where the walls of the two atria are in contact. The figure displays the result of interpolating, from the coarse to the fine mesh, a displacement field that maps the starting configuration onto the stress-free reference configuration, computed during preprocessing as in \cite{fedele2023comprehensive,barnafi2024reconstructing}.

We observe that in both regions the \ac{RLRBF} interpolation yields significant oscillations. \rlrbfg{} without high-curvature detection (i.e. with $\beta = +\infty$) eliminates the oscillations near the left atrial appendage (top row of \cref{fig:high-curvature}), but not those near the fossa ovalis (bottom row of \cref{fig:high-curvature}). Conversely, introducing high-curvature detection by setting $\beta = \num{0.5}$ removes spurious oscillations in both regions.

An intuitive explanation for this behavior is that the condition $\beta h_{\max} + \|\mathbf x - \mathbf y\| < \discgeodist(\mathbf x, \mathbf y)$ in \eqref{eq:distance-cutoff} empirically distinguishes between cases in which the exact geodesic distance $\geodist$ is better approximated by the Euclidean distance (if the condition is not met) or by discrete geodesic distance (if the condition is met). In the latter case, it is preferable to evaluate $d(\mathbf x, \mathbf y; R)$ through $\discgeodist(\mathbf x, \mathbf y)$ (instead of $\|\mathbf x - \mathbf y\|$), despite $\discgeodist$ being a potentially inaccurate approximation of $\geodist$.

\section{Conclusions}
\label{sec:conclusions}

We proposed an interpolation method based on \ac{RLRBF} that uses a discrete approximation of the geodesic distance to take into account the geometry of the interpolation domain. The interpolation method, referred to as \rlrbfg{}, allows to transfer variables between non-matching meshes, and prevents the onset of spurious oscillations associated with holes, cuts and otherwise complex geometrical features. Our numerical experiments demonstrate that the proposed method overcomes limitations of \ac{RLRBF}, while retaining good convergence properties.

Furthermore, we focused on the distributed-memory parallel implementation of \rlrbfg{} interpolation, a non-trivial task due to the two meshes having independent parallel partitionings. The scalability tests carried out in this work indicate that the proposed implementation, available as open source within \lifex{} \cite{africa2022lifexcore}, has nearly ideal scalability properties, and is therefore well suited to large scale applications.

Finally, we demonstrate the effectiveness of \rlrbfg{} interpolation in a realistic application to a whole-heart cardiac electromechanics model. This application is a particularly challenging setting for the proposed method, given the geometrical complexity of the domain. This test shows that the \rlrbfg{} method is capable of transferring data between arbitrarily refined and even non-nested meshes, with an acceptable impact on the computational costs, thereby providing an effective tool in developing multiphysics computational models in geometrically complex domains.

\appendix
\section{Whole-heart electromechanics model}
\label{app:em-model}

We briefly report the equations of the electromechanical model used in the numerical experiments of \cref{sec:heart}. The model used in this work was introduced in \cite{fedele2023comprehensive}, and we refer to the original work for additional details.

The domain $\Omega$, representing a four-chamber realistic human heart \cite{zygote}, is decomposed into several subsets: $\Omega\la$ (left atrium),  $\Omega\ra$ (right atrium), $\Omega\ventricles$ (ventricles), $\Omega\ao$ (ascending aorta), $\Omega\pt$ (pulmonary trunk) and $\Omega_\text{caps}$ (solid caps on vein inlets and valves). We define the myocardial subdomain as $\Omega\myo = \Omega\la \cup \Omega\ra \cup \Omega\ventricles$.

The domain boundary $\partial\Omega$ is also decomposed into several subsets: $\Gamma_\text{rings}$ (the artificial boundaries where the veins and arteries are cut), $\Gamma_\text{epi}^\text{PF}$ (the epicardial portion of the boundary in contact with the pericardial fluid), $\Gamma_\text{epi}^\text{EAT}$ (the epicardial portion of the boundary in contact with the epicardial adipous tissue), $\Gamma_\text{epi}^\text{art}$ (the epithelial surface of the two arteries), the endocardial surfaces of the four chambers, denoted by $\Gamma\la$, $\Gamma\lv$, $\Gamma\ra$, $\Gamma\rv$, and the endothelial surfaces of the two arteries, $\Gamma\ao$ and $\Gamma\pt$.

The domain $\Omega$ is obtained starting from a known configuration corresponding to medical imaging data. We consider the Zygote Heart Model \cite{zygote}, preprocessed with VMTK using the algorithms described in \cite{fedele2021polygonal}. The model represents a configuration in which the heart muscle is subject to the pressure of the blood it contains and of the organs that surround it. To recover the stress-free configuration $\Omega$, we use a reference configuration recovery algorithm described in \cite{regazzoni2022cardiac, barnafi2024reconstructing}.

The model has the following unknowns:
\begin{align*}
  \potential &: \Omega\myo \times(0, T) \to \mathbb{R} & \text{transmembrane potential}\;, \\
  \ionicvars\ttp &: \Omega\ventricles\times(0, T) \to \mathbb{R}^{N_\text{TTP}} & \text{state variables of the \ac{TTP} model}\;, \\
  \ionicvars\crn &: \left(\Omega\la \cup \Omega\ra\right)\times(0, T) \to \mathbb{R}^{N_\text{CRN}} & \text{state variables of the \ac{CRN} model}\;, \\
  \actvars &: \domain\myo\times(0, T) \to \mathbb{R}^{N_\text{act}} & \text{activation variables}\;, \\
  \displacement &: \domain\times(0, T) \to \mathbb{R}^3 & \text{solid displacement}\;, \\
  \circvars &: (0, T) \to \mathbb{R}^{N_\text{circ}} & \text{circulation state variables}\;.
\end{align*}
We reconstruct the orientation of cardiac muscle fibers through a \ac{LDRBM} \cite{piersanti2021modeling, fedele2023comprehensive}, defining an orthonormal triplet $\{\fibers, \sheets, \normals\}$ at every point of $\Omega$.

The evolution of the ionic variables is prescribed by the \ac{TTP} and \ac{CRN} ionic models.

In the myocardial subdomain $\Omega\myo$ we solve the monodomain equation with geometry-mediated mechano-electric feedback effects \cite{collifranzone2014mathematical, salvador2022role}:
\begin{equation}
  \begin{cases}
    \jacobian \chi C_\text{m} \pdv{\potential}{t}
      - \div(\jacobian \Ftensor^{-1}\difftensor\Ftensor^{-T} \grad\potential)
      + \jacobian\chi\Iion(\potential, \ionicvars\ttp, \ionicvars\crn)
      = \jacobian\chi\Iapp & \text{in } \Omega\myo \times (0,T)\;, \\
    \jacobian\Ftensor^{-1}\difftensor\Ftensor^{-T}\grad\potential\cdot\normal = 0 & \text{on }\partial\Omega\myo \times (0, T)\;, \\
    \potential = \potential_0 & \text{in }\Omega\myo \times \{0\}\;,
  \end{cases}
  \label{eq:monodomain}
\end{equation}
wherein $\chi$ is the membrane area-to-surface ratio, $C_\text{m}$ is the membrane capacitance, and $\difftensor$ is the conductivity tensor, defined through
\begin{equation}
  \difftensor =
  \sigmaf\frac{\Ftensor\fibers\otimes\Ftensor\fibers}{\|\Ftensor\fibers\|^2} +
  \sigmas\frac{\Ftensor\sheets\otimes\Ftensor\sheets}{\|\Ftensor\sheets\|^2} +
  \sigman\frac{\Ftensor\normals\otimes\Ftensor\normals}{\|\Ftensor\normals\|^2}\;,
\end{equation}
with $\sigmaf$, $\sigmas$ and $\sigman$ positive coefficient with different constant values in $\Omega\ventricles$, $\Omega\la$ and $\Omega\ra$.

The ionic current term $\Iion(\potential, \ionicvars\ttp, \ionicvars\crn)$ is defined according to a different ionic model, depending on the subdomain:
\begin{equation}
  \Iion(\potential, \ionicvars\ttp, \ionicvars\crn) =
  \begin{cases}
    \Iion^\text{TTP}(\potential, \ionicvars\ttp) & \text{in } \Omega\ventricles\;, \\
    \Iion^\text{CRN}(\potential, \ionicvars\crn) & \text{in } \Omega\la \cup \Omega\ra\;.
  \end{cases}
\end{equation}
The evolution of the ionic variables $\ionicvars\ttp$ and $\ionicvars\crn$ is prescribed by a system of \acp{ODE} at every point in the corresponding subdomain:
\begin{align}
  & \begin{cases}
    \pdv{\ionicvars\ttp}{t} = \mathbf{F}_\text{ion}^\text{TTP}(v, \ionicvars\ttp) & \text{in } \Omega\ventricles \times (0, T)\;, \\
    \ionicvars\ttp = \ionicvars_{\text{TTP},0} & \text{in } \Omega\ventricles \times \{0\}\;,
  \end{cases} \\
  & \begin{cases}
    \pdv{\ionicvars\crn}{t} = \mathbf{F}_\text{ion}^\text{CRN}(v, \ionicvars\crn) & \text{in } \left(\Omega\la \cup \Omega\ra\right) \times (0, T)\;, \\
    \ionicvars\crn = \ionicvars_{\text{CRN},0} & \text{in } \left(\Omega\la \cup \Omega\ra\right) \times \{0\}\;.
  \end{cases}
\end{align}
The definitions of $\Iion^\text{TTP}$, $\Iion^\text{CRN}$, $\mathbf F_\text{ion}^\text{TTP}$ and $\mathbf F_\text{ion}^\text{CRN}$ are given in the original works \cite{ten2006alternans, courtemanche1998ionic}. The initial conditions for $v$, $\ionicvars\ttp$ and $\ionicvars\crn$ are determined as the limit cycle of a zero-dimensional simulation, as in \cite{fedele2023comprehensive, piersanti20223d}.

Both ionic models include a variable representing the intracellular calcium concentration, denoted by $w\ttp^{\text{Ca}}$ and $w\crn^{\text{Ca}}$. For convenience, we define
\begin{equation}
  \calcium = \begin{cases}
    w\ttp^{\text{Ca}} & \text{in } \Omega\ventricles\;, \\
    w\crn^{\text{Ca}} & \text{in } \Omega\la \cup \Omega\ra\;. \\
  \end{cases}
\end{equation}
The generation of contractile force is modeled through the RDQ20-MF model \cite{regazzoni2020biophysically}, which can be expressed as a system of \acp{ODE}:
\begin{equation}
  \begin{cases}
    \pdv{\actvars}{t} = \mathbf{F}_\text{act}\left(\actvars, \calcium, \SL, \pdv{\SL}{t}\right) & \text{in } \Omega\myo \times (0, T)\;, \\
    \actvars = \actvars_0 & \text{in } \Omega\myo \times \{0\}\;,
  \end{cases}
\end{equation}
where $\SL$ denotes the sarcomere length, defined as $\SL = \SL_0\|\Ftensor\fibers\|$. The force generation model also defines the active stress tensor $\actstress = \actstress(\actvars, \displacement)$. The initial condition for $\actvars$ is the steady state of a zero-dimensional simulation \cite{fedele2023comprehensive, piersanti20223d}.

The deformation of the muscle is governed by the elastodynamics equation, endowed with suitable boundary conditions \cite{fedele2023comprehensive}:
\begin{equation}
  \begin{cases}
    \pdv[2]{\displacement}{t} - \div\stress(\actvars, \displacement) = \mathbf 0 & \text{in } \Omega \times (0, T)\;, \\
    \stress(\actvars, \displacement)\normal + (\normal \otimes \normal)\left(\Kepi \displacement + \Cepi \pdv{\displacement}{t}\right) = \mathbf 0 & \text{on } \left(\Gamma_\text{epi}^{PF} \cup \Gamma_\text{epi}^\text{EAT} \cup \Gamma_\text{epi}^\text{art}\right) \times (0, T)\;,\\
    \stress(\actvars, \displacement)\normal = -p_\text{i}(t)\normal & \text{on } \Gamma_\text{i} \times (0, T) \qquad \text{for } \text{i} \in \{\text{LA}, \text{LV}, \text{RA}, \text{RV}, \text{Ao}, \text{PT}\}\;, \\
    \displacement = \mathbf{0} & \text{on } \Gamma_\text{rings} \times (0, T)\;, \\
    \displacement = \displacement_0 & \text{in } \Omega \times \{ 0 \}\;, \\
    \pdv{\displacement}{t} = \mathbf{0} & \text{in } \Omega \times \{ 0 \}\;.
  \end{cases}
\end{equation}
The initial displacement is determined through the solution of a quasi-static problem imposing suitable initial endocardial pressures \cite{regazzoni2022cardiac, fedele2023comprehensive}.

Finally, the evolution of the endocardial pressures $p_\text{i}$ is regulated by a lumped-parameter model of the circulatory system \cite{blanco20103d, hirschvogel2017monolithic, regazzoni2022cardiac, piersanti20223d, fedele2023comprehensive}, expressed by a system of ODEs:
\begin{equation}
  \begin{cases}
    \dv{\circvars}{t} = \mathbf{F}_\text{circ}(\circvars, t) & \text{in } (0, T)\;, \\
    \circvars(0) = \circvars_0\;.
  \end{cases}
\end{equation}
The circulation and mechanics models are further coupled by imposing the volume constraints
\begin{equation}
  V^\text{3D}_\text{i}(t) = V^\text{0D}_\text{i}(t) \quad \text{for } \text{i} \in \{\text{LA}, \text{LV}, \text{RA}, \text{RV}\}\;,
\end{equation}
where $V^\text{3D}_\text{i}$ and $V^\text{0D}_\text{i}$ are the volumes of chamber $\text{i}$ computed through the 3D mechanics model and the circulation model respectively \cite{fedele2023comprehensive}.

\section*{Acknowledgements}

M.B. acknowledges W. Bangerth (Colorado State University, Fort Collins, CO, USA) and M. Feder (SISSA, Trieste, Italy) for their helpful suggestions on parallel geometric algorithms and data structures.

The present research is part of the activities of ``Dipartimento di Eccellenza
2023--2027'', MUR, Italy, Dipartimento di Matematica, Politecnico di Milano. M.B. and L.D. have received support from the project PRIN2022, MUR, Italy, 2023--2025, 202232A8AN ``Computational modeling of the heart: from efficient numerical solvers to cardiac digital twins''. F.R. has received support from the project PRIN2022, MUR, Italy, 2023--2025, P2022N5ZNP ``SIDDMs: shape-informed data-driven models for parametrized PDEs, with application to computational cardiology''. The authors acknowledge their membership to INdAM GNCS - Gruppo Nazionale per il Calcolo Scientifico (National Group for Scientific Computing, Italy). This project has been partially supported by the INdAM-GNCS Project CUP\textunderscore E53C22001930001.

The authors acknowledge the CINECA award under the ISCRA initiative, for the availability of high performance computing resources and support.

\bibliographystyle{abbrv}
\bibliography{bibliography}

\begin{thebibliography}{10}

\bibitem{africa2022lifexcore}
P.~C. Africa.
\newblock lifex: A flexible, high performance library for the numerical
  solution of complex finite element problems.
\newblock {\em SoftwareX}, 20:101252, 2022.

\bibitem{africa2024lifexcfd}
P.~C. Africa, I.~Fumagalli, M.~Bucelli, A.~Zingaro, M.~Fedele, L.~Dede', and
  A.~Quarteroni.
\newblock lifex-cfd: An open-source computational fluid dynamics solver for
  cardiovascular applications.
\newblock {\em Computer Physics Communications}, 296:109039, 2024.

\bibitem{africa2023lifexfiber}
P.~C. Africa, R.~Piersanti, M.~Fedele, L.~Dede', and A.~Quarteroni.
\newblock lifex-fiber: an open tool for myofibers generation in cardiac
  computational models.
\newblock {\em BMC Bioinformatics}, 24:143, 2023.

\bibitem{africa2023lifexep}
P.~C. Africa, R.~Piersanti, F.~Regazzoni, M.~Bucelli, M.~Salvador, M.~Fedele,
  S.~Pagani, L.~Dede', and A.~Quarteroni.
\newblock lifex-ep: a robust and efficient software for cardiac
  electrophysiology simulations.
\newblock {\em BMC bioinformatics}, 24(1):389, 2023.

\bibitem{arndt2020dealii}
D.~Arndt, W.~Bangerth, D.~Davydov, T.~Heister, L.~Heltai, M.~Kronbichler,
  M.~Maier, J.~Pelteret, B.~Turcksin, and D.~Wells.
\newblock The \texttt{deal.II} finite element library: design, features, and
  insights.
\newblock {\em Computers \& Mathematics with Applications}, 81:407--422, 2021.

\bibitem{arndt2022dealii9.4}
D.~Arndt, W.~Bangerth, M.~Feder, M.~Fehling, R.~Gassm{\"o}ller, T.~Heister,
  L.~Heltai, M.~Kronbichler, M.~Maier, P.~Munch, J.-P. Pelteret, S.~Sticko,
  B.~Turcksin, and D.~Wells.
\newblock The \texttt{deal.II} library, version 9.4.
\newblock {\em Journal of Numerical Mathematics}, 30(3):231--246, 2022.

\bibitem{augustin2016patient}
C.~M. Augustin, A.~Crozier, A.~Neic, A.~J. Prassl, E.~Karabelas, T.~Ferreira~da
  Silva, J.~F. Fernandes, F.~Campos, T.~Kuehne, and G.~Plank.
\newblock Patient-specific modeling of left ventricular electromechanics as a
  driver for haemodynamic analysis.
\newblock {\em EP Europace}, 18:iv121--iv129, 2016.

\bibitem{augustin2021computationally}
C.~M. Augustin, M.~A. Gsell, E.~Karabelas, E.~Willemen, F.~W. Prinzen,
  J.~Lumens, E.~J. Vigmond, and G.~Plank.
\newblock A computationally efficient physiologically comprehensive {3D--0D}
  closed-loop model of the heart and circulation.
\newblock {\em Computer methods in applied mechanics and engineering},
  386:114092, 2021.

\bibitem{barnafi2024reconstructing}
N.~Barnafi, F.~Regazzoni, and D.~Riccobelli.
\newblock Reconstructing relaxed configurations in elastic bodies: Mathematical
  formulations and numerical methods for cardiac modeling.
\newblock {\em Computer Methods in Applied Mechanics and Engineering},
  423:116845, 2024.

\bibitem{beatson1999fast}
R.~K. Beatson, J.~B. Cherrie, and C.~T. Mouat.
\newblock Fast fitting of radial basis functions: Methods based on
  preconditioned {GMRES} iteration.
\newblock {\em Advances in Computational Mathematics}, 11:253--270, 1999.

\bibitem{blanco20103d}
P.~J. Blanco and R.~A. Feij{\'o}o.
\newblock A {3D-1D-0D} computational model for the entire cardiovascular
  system.
\newblock {\em Mec{\'a}nica Computacional}, 29(59):5887--5911, 2010.

\bibitem{boost}
Official boost website.
\newblock \url{https://www.boost.org/} (last accessed: 6 November 2023).

\bibitem{brown2005approximate}
D.~Brown, L.~Ling, E.~Kansa, and J.~Levesley.
\newblock On approximate cardinal preconditioning methods for solving {PDEs}
  with radial basis functions.
\newblock {\em Engineering Analysis with Boundary Elements}, 29(4):343--353,
  2005.

\bibitem{bucelli2023preserving}
M.~Bucelli, F.~Regazzoni, L.~Dede', and A.~Quarteroni.
\newblock Preserving the positivity of the deformation gradient determinant in
  intergrid interpolation by combining {RBFs} and {SVD}: application to cardiac
  electromechanics.
\newblock {\em Computer Methods in Applied Mechanics and Engineering}, page
  116292, 2023.

\bibitem{chourdakis2022precice}
G.~Chourdakis, K.~Davis, B.~Rodenberg, M.~Schulte, F.~Simonis, B.~Uekermann,
  G.~Abrams, H.-J. Bungartz, L.~C. Yau, I.~Desai, et~al.
\newblock {preCICE} v2: A sustainable and user-friendly coupling library.
\newblock {\em Open Research Europe}, 2, 2022.

\bibitem{ciaccia1997indexing}
P.~Ciaccia, M.~Patella, F.~Rabitti, and P.~Zezula.
\newblock Indexing metric spaces with m-tree.
\newblock In {\em SEBD}, volume~97, pages 67--86, 1997.

\bibitem{clevenger2020flexible}
T.~C. Clevenger, T.~Heister, G.~Kanschat, and M.~Kronbichler.
\newblock A flexible, parallel, adaptive geometric multigrid method for {FEM}.
\newblock {\em ACM Transactions on Mathematical Software (TOMS)}, 47(1):1--27,
  2020.

\bibitem{collifranzone2014mathematical}
P.~Colli~Franzone, L.~F. Pavarino, and S.~Scacchi.
\newblock {\em {Mathematical Cardiac Electrophysiology}}, volume~13.
\newblock Springer, 2014.

\bibitem{courtemanche1998ionic}
M.~Courtemanche, R.~J. Ramirez, and S.~Nattel.
\newblock Ionic mechanisms underlying human atrial action potential properties:
  insights from a mathematical model.
\newblock {\em American Journal of Physiology-Heart and Circulatory
  Physiology}, 275(1):H301--H321, 1998.

\bibitem{demarchi2020convergence}
S.~De~Marchi and H.~Wendland.
\newblock On the convergence of the rescaled localized radial basis function
  method.
\newblock {\em Applied Mathematics Letters}, 99:105996, 2020.

\bibitem{dealii}
Official \texttt{deal.ii} website.
\newblock \url{https://www.dealii.org/} (last accessed: 6 November 2023).

\bibitem{deparis2016internodes}
S.~Deparis, D.~Forti, P.~Gervasio, and A.~Quarteroni.
\newblock {INTERNODES}: an accurate interpolation-based method for coupling the
  {Galerkin} solutions of {PDEs} on subdomains featuring non-conforming
  interfaces.
\newblock {\em Computers \& Fluids}, 141:22--41, 2016.

\bibitem{deparis2014rescaled}
S.~Deparis, D.~Forti, and A.~Quarteroni.
\newblock A rescaled localized radial basis function interpolation on
  non-cartesian and nonconforming grids.
\newblock {\em SIAM Journal on Scientific Computing}, 36(6):A2745--A2762, 2014.

\bibitem{dijkstra1959note}
E.~W. Dijkstra.
\newblock A note on two problems in connexion with graphs.
\newblock {\em Numerische Mathematik}, 1(1):269--271, 1959.

\bibitem{farah2016volumetric}
P.~Farah, A.-T. Vuong, W.~Wall, and A.~Popp.
\newblock Volumetric coupling approaches for multiphysics simulations on
  non-matching meshes.
\newblock {\em International Journal for Numerical Methods in Engineering},
  108(12):1550--1576, 2016.

\bibitem{fedele2023comprehensive}
M.~Fedele, R.~Piersanti, F.~Regazzoni, M.~Salvador, P.~C. Africa, M.~Bucelli,
  A.~Zingaro, L.~Dede', and A.~Quarteroni.
\newblock A comprehensive and biophysically detailed computational model of the
  whole human heart electromechanics.
\newblock {\em Computer Methods in Applied Mechanics and Engineering},
  410:115983, 2023.

\bibitem{fedele2021polygonal}
M.~Fedele and A.~Quarteroni.
\newblock Polygonal surface processing and mesh generation tools for the
  numerical simulation of the cardiac function.
\newblock {\em International Journal for Numerical Methods in Biomedical
  Engineering}, 37(4):e3435, 2021.

\bibitem{franke1998solving}
C.~Franke and R.~Schaback.
\newblock Solving partial differential equations by collocation using radial
  basis functions.
\newblock {\em Applied Mathematics and Computation}, 93(1):73--82, 1998.

\bibitem{gerach2024differential}
T.~Gerach and A.~Loewe.
\newblock Differential effects of mechano-electric feedback mechanisms on
  whole-heart activation, repolarization, and tension.
\newblock {\em The Journal of Physiology}, 2024.

\bibitem{gerach2021electro}
T.~Gerach, S.~Schuler, J.~Fr{\"o}hlich, L.~Lindner, E.~Kovacheva, R.~Moss,
  E.~M. W{\"u}lfers, G.~Seemann, C.~Wieners, and A.~Loewe.
\newblock Electro-mechanical whole-heart digital twins: a fully coupled
  multi-physics approach.
\newblock {\em Mathematics}, 9(11):1247, 2021.

\bibitem{gillette2021framework}
K.~Gillette, M.~A. Gsell, A.~J. Prassl, E.~Karabelas, U.~Reiter, G.~Reiter,
  T.~Grandits, C.~Payer, D.~{\v{S}}tern, M.~Urschler, J.~D. Bayer, C.~M.
  Augustin, A.~Neic, T.~Pock, E.~J. Vigmond, and G.~Plank.
\newblock A framework for the generation of digital twins of cardiac
  electrophysiology from clinical 12-leads {ECGs}.
\newblock {\em Medical Image Analysis}, 71:102080, 2021.

\bibitem{gurev2011models}
V.~Gurev, T.~Lee, J.~Constantino, H.~Arevalo, and N.~A. Trayanova.
\newblock Models of cardiac electromechanics based on individual hearts imaging
  data: Image-based electromechanical models of the heart.
\newblock {\em Biomechanics and Modeling in Mechanobiology}, 10(3):295--306,
  2011.

\bibitem{guttman1984rtrees}
A.~Guttman.
\newblock R-trees: A dynamic index structure for spatial searching.
\newblock In {\em Proceedings of the 1984 ACM SIGMOD international conference
  on Management of data}, pages 47--57, 1984.

\bibitem{heijman2021computational}
J.~Heijman, H.~Sutanto, H.~J. Crijns, S.~Nattel, and N.~A. Trayanova.
\newblock Computational models of atrial fibrillation: Achievements,
  challenges, and perspectives for improving clinical care.
\newblock {\em Cardiovascular Research}, 117(7):1682--1699, 2021.

\bibitem{hirschvogel2017monolithic}
M.~Hirschvogel, M.~Bassilious, L.~Jagschies, S.~M. Wildhirt, and M.~W. Gee.
\newblock A monolithic {3D-0D} coupled closed-loop model of the heart and the
  vascular system: experiment-based parameter estimation for patient-specific
  cardiac mechanics.
\newblock {\em International Journal for Numerical Methods in Biomedical
  Engineering}, 33(8):e2842, 2017.

\bibitem{karabelas2022accurate}
E.~Karabelas, M.~A. Gsell, G.~Haase, G.~Plank, and C.~M. Augustin.
\newblock An accurate, robust, and efficient finite element framework with
  applications to anisotropic, nearly and fully incompressible elasticity.
\newblock {\em Computer Methods in Applied Mechanics and Engineering},
  394:114887, 2022.

\bibitem{kronbichler2019multigrid}
M.~Kronbichler and K.~Ljungkvist.
\newblock Multigrid for matrix-free high-order finite element computations on
  graphics processors.
\newblock {\em ACM Transactions on Parallel Computing (TOPC)}, 6(1):1--32,
  2019.

\bibitem{lanthier2001approximating}
M.~Lanthier, A.~Maheshwari, and J.~R.~Sack.
\newblock Approximating shortest paths on weighted polyhedral surfaces.
\newblock {\em Algorithmica}, 30:527--562, 2001.

\bibitem{lanthier1997approximating}
M.~Lanthier, A.~Maheshwari, and J.-R. Sack.
\newblock Approximating weighted shortest paths on polyhedral surfaces.
\newblock In {\em Proceedings of the thirteenth annual symposium on
  Computational geometry}, pages 274--283, 1997.

\bibitem{levrero2020sensitivity}
F.~Levrero-Florencio, F.~Margara, E.~Zacur, A.~Bueno-Orovio, Z.~Wang,
  A.~Santiago, J.~Aguado-Sierra, G.~Houzeaux, V.~Grau, D.~Kay, M.~V\'{a}zquez,
  R.~Ruiz-Baier, and B.~Rodriguez.
\newblock Sensitivity analysis of a strongly-coupled human-based
  electromechanical cardiac model: Effect of mechanical parameters on
  physiologically relevant biomarkers.
\newblock {\em Computer Methods in Applied Mechanics and Engineering},
  361:112762, 2020.

\bibitem{lifex}
Official \texttt{life\textsuperscript{x}} website.
\newblock \url{https://lifex.gitlab.io/} (last accessed: 6 November 2023).

\bibitem{manolopoulos2006rtrees}
Y.~Manolopoulos, A.~N. Papadopoulos, and Y.~Theodoridis.
\newblock {\em R-Trees: Theory and Applications}.
\newblock Springer Science \& Business Media, 2006.

\bibitem{morse2005interpolating}
B.~S. Morse, T.~S. Yoo, P.~Rheingans, D.~T. Chen, and K.~R. Subramanian.
\newblock Interpolating implicit surfaces from scattered surface data using
  compactly supported radial basis functions.
\newblock In {\em ACM SIGGRAPH 2005 Courses}, pages 78--es. 2005.

\bibitem{pagani2021computational}
S.~Pagani, L.~Dede', A.~Frontera, M.~Salvador, L.~Limite, A.~Manzoni,
  F.~Lipartiti, G.~Tsitsinakis, A.~Hadjis, P.~Della~Bella, and A.~Quarteroni.
\newblock A computational study of the electrophysiological substrate in
  patients suffering from atrial fibrillation.
\newblock {\em Frontiers in Physiology}, 12:673612, 2021.

\bibitem{peirlinck2021precision}
M.~Peirlinck, F.~S. Costabal, J.~Yao, J.~Guccione, S.~Tripathy, Y.~Wang,
  D.~Ozturk, P.~Segars, T.~Morrison, S.~Levine, and E.~Kuhl.
\newblock Precision medicine in human heart modeling.
\newblock {\em Biomechanics and Modeling in Mechanobiology}, 20(3):803--831,
  2021.

\bibitem{piersanti2021modeling}
R.~Piersanti, P.~C. Africa, M.~Fedele, C.~Vergara, L.~Dede', A.~F. Corno, and
  A.~Quarteroni.
\newblock Modeling cardiac muscle fibers in ventricular and atrial
  electrophysiology simulations.
\newblock {\em Computer Methods in Applied Mechanics and Engineering},
  373:113468, 2021.

\bibitem{piersanti20223d}
R.~Piersanti, F.~Regazzoni, M.~Salvador, A.~F. Corno, C.~Vergara, and
  A.~Quarteroni.
\newblock {3D--0D} closed-loop model for the simulation of cardiac
  biventricular electromechanics.
\newblock {\em Computer Methods in Applied Mechanics and Engineering},
  391:114607, 2022.

\bibitem{quarteroni2019mathematical}
A.~Quarteroni, L.~Dede', A.~Manzoni, and C.~Vergara.
\newblock {\em Mathematical modelling of the human cardiovascular system: data,
  numerical approximation, clinical applications}, volume~33.
\newblock Cambridge University Press, 2019.

\bibitem{quarteroni2022modeling}
A.~Quarteroni, L.~Dede', and F.~Regazzoni.
\newblock Modeling the cardiac electromechanical function: A mathematical
  journey.
\newblock {\em Bulletin of the American Mathematical Society}, 59(3):371--403,
  2022.

\bibitem{quarteroni2023mathematical}
A.~Quarteroni, L.~Dede', F.~Regazzoni, and C.~Vergara.
\newblock A mathematical model of the human heart suitable to address clinical
  problems.
\newblock {\em Japan Journal of Industrial and Applied Mathematics}, pages
  1--21, 2023.

\bibitem{regazzoni2020biophysically}
F.~Regazzoni, L.~Dede', and A.~Quarteroni.
\newblock Biophysically detailed mathematical models of multiscale cardiac
  active mechanics.
\newblock {\em PLoS Computational Biology}, 16(10):e1008294, 2020.

\bibitem{regazzoni2022cardiac}
F.~Regazzoni, M.~Salvador, P.~C. Africa, M.~Fedele, L.~Dede', and
  A.~Quarteroni.
\newblock A cardiac electromechanical model coupled with a lumped-parameter
  model for closed-loop blood circulation.
\newblock {\em Journal of Computational Physics}, 457:111083, 2022.

\bibitem{saad2003iterative}
Y.~Saad.
\newblock {\em {Iterative Methods for Sparse Linear Systems}}.
\newblock SIAM, 2003.

\bibitem{sala2005robust}
M.~Sala and M.~A. Heroux.
\newblock Robust algebraic preconditioners using {IFPACK} 3.0.
\newblock Technical report, Sandia National Lab.(SNL-NM), Albuquerque, NM
  (United States), 2005.

\bibitem{salvador2020intergrid}
M.~Salvador, L.~Dede', and A.~Quarteroni.
\newblock An intergrid transfer operator using radial basis functions with
  application to cardiac electromechanics.
\newblock {\em Computational Mechanics}, 66:491--511, 2020.

\bibitem{salvador2021electromechanical}
M.~Salvador, M.~Fedele, P.~C. Africa, E.~Sung, L.~Dede', A.~Prakosa,
  J.~Chrispin, N.~A. Trayanova, and A.~Quarteroni.
\newblock Electromechanical modeling of human ventricles with ischemic
  cardiomyopathy: numerical simulations in sinus rhythm and under arrhythmia.
\newblock {\em Computers in Biology and Medicine}, 136:104674, 2021.

\bibitem{salvador2022role}
M.~Salvador, F.~Regazzoni, S.~Pagani, N.~A. Trayanova, and A.~Quarteroni.
\newblock The role of mechano-electric feedbacks and hemodynamic coupling in
  scar-related ventricular tachycardia.
\newblock {\em Computers in Biology and Medicine}, 142:105203, 2022.

\bibitem{stella2022fast}
S.~Stella, F.~Regazzoni, C.~Vergara, L.~Dede', and A.~Quarteroni.
\newblock A fast cardiac electromechanics model coupling the eikonal and the
  nonlinear mechanics equations.
\newblock {\em Mathematical Models and Methods in Applied Sciences},
  32(08):1531--1556, 2022.

\bibitem{strocchi2020simulating}
M.~Strocchi, M.~A. Gsell, C.~M. Augustin, O.~Razeghi, C.~H. Roney, A.~J.
  Prassl, E.~J. Vigmond, J.~M. Behar, J.~S. Gould, C.~A. Rinaldi, M.~J. Bishop,
  G.~Plank, and S.~A. Niederer.
\newblock Simulating ventricular systolic motion in a four-chamber heart model
  with spatially varying {Robin} boundary conditions to model the effect of the
  pericardium.
\newblock {\em Journal of Biomechanics}, 101:109645, 2020.

\bibitem{sundnes2007computing}
J.~Sundnes, G.~T. Lines, X.~Cai, B.~F. Nielsen, K.-A. Mardal, and A.~Tveito.
\newblock {\em {Computing the Electrical Activity in the Heart}}, volume~1.
\newblock Springer Science \& Business Media, 2007.

\bibitem{sung2020personalized}
E.~Sung, A.~Prakosa, K.~N. Aronis, S.~Zhou, S.~L. Zimmerman, H.~Tandri,
  S.~Nazarian, R.~D. Berger, J.~Chrispin, and N.~A. Trayanova.
\newblock Personalized digital-heart technology for ventricular tachycardia
  ablation targeting in hearts with infiltrating adiposity.
\newblock {\em Circulation: Arrhythmia and Electrophysiology}, 13(12):e008912,
  2020.

\bibitem{ten2006alternans}
K.~H. Ten~Tusscher and A.~V. Panfilov.
\newblock Alternans and spiral breakup in a human ventricular tissue model.
\newblock {\em American Journal of Physiology-Heart and Circulatory
  Physiology}, 291(3):H1088--H1100, 2006.

\bibitem{viola2023gpu}
F.~Viola, G.~Del~Corso, R.~De~Paulis, and R.~Verzicco.
\newblock Gpu accelerated digital twins of the human heart open new routes for
  cardiovascular research.
\newblock {\em Scientific Reports}, 13(1):8230, 2023.

\bibitem{voet2022internodes}
Y.~Voet, G.~Anciaux, S.~Deparis, and P.~Gervasio.
\newblock The {INTERNODES} method for applications in contact mechanics and
  dedicated preconditioning techniques.
\newblock {\em Computers \& Mathematics with Applications}, 127:48--64, 2022.

\bibitem{washio2013multiscale}
T.~Washio, J.-i. Okada, A.~Takahashi, K.~Yoneda, Y.~Kadooka, S.~Sugiura, and
  T.~Hisada.
\newblock Multiscale heart simulation with cooperative stochastic cross-bridge
  dynamics and cellular structures.
\newblock {\em Multiscale Modeling \& Simulation}, 11(4):965--999, 2013.

\bibitem{wendland1995piecewise}
H.~Wendland.
\newblock Piecewise polynomial, positive definite and compactly supported
  radial functions of minimal degree.
\newblock {\em Advances in computational Mathematics}, 4(1):389--396, 1995.

\bibitem{wendland1999meshless}
H.~Wendland.
\newblock Meshless {Galerkin} methods using radial basis functions.
\newblock {\em Mathematics of computation}, 68(228):1521--1531, 1999.

\bibitem{williams1964algorithm}
J.~W.~J. Williams.
\newblock Algorithm 232 - heapsort.
\newblock {\em Communications of the ACM}, 7(6):347--348, 1964.

\bibitem{zappon2022reduced}
E.~Zappon, A.~Manzoni, P.~Gervasio, and A.~Quarteroni.
\newblock A reduced order model for domain decompositions with non-conforming
  interfaces.
\newblock {\em arXiv preprint arXiv:2206.09618}, 2022.

\bibitem{zappon2023efficient}
E.~Zappon, A.~Manzoni, and A.~Quarteroni.
\newblock Efficient and certified solution of parametrized one-way coupled
  problems through {DEIM}-based data projection across non-conforming
  interfaces.
\newblock {\em Advances in Computational Mathematics}, 49(2):21, 2023.

\bibitem{zingaro2023electromechanics}
A.~Zingaro, M.~Bucelli, R.~Piersanti, F.~Regazzoni, L.~Dede', and
  A.~Quarteroni.
\newblock An electromechanics-driven fluid dynamics model for the simulation of
  the whole human heart.
\newblock {\em Journal of Computational Physics (in press)}, 2024.

\bibitem{zygote}
{Zygote Media Group Inc.}
\newblock Zygote solid {3D} heart generation {II} developement report.
  {Technical report.}
\newblock 2014.

\end{thebibliography}

\end{document}